\documentclass{article}

\usepackage{color}
\usepackage{geometry}
\usepackage{amsmath}
\usepackage{amssymb}
\usepackage{amsthm}
\usepackage[all]{xy}
\usepackage{tikz}
\usetikzlibrary{calc}

\newtheorem{theorem}{Theorem}[section]
\newtheorem{prop}[theorem]{Proposition}
\newtheorem{lemma}[theorem]{Lemma}
\newtheorem{defi}[theorem]{Definition}
\newtheorem{example}[theorem]{Example}
\newtheorem{rem}[theorem]{Remark}

\newtheorem{cor}[theorem]{Corollary}

\newcommand{\ord}{\mathord{\mathrm{ord}}}

\newcommand{\PP}{\mathbb{P}}
\newcommand{\CC}{\mathbb{C}}

\newcommand{\ZZ}{\mathbb{Z}}
\newcommand{\RR}{\mathbb{R}}

\title{Projectivized ASEP-Bethe system and counting solutions} 
\author{Shinsuke Iwao\\[1em]
\normalsize
Faculty of Business and Commerce, 
Keio University \\
\normalsize
4-1-1 Hiyoshi, Kohoku-ku, Yokohama, Kanagawa 223-8521, Japan\\
\normalsize
Email: iwao-s@keio.jp
}
\date{\today}
\begin{document}

\maketitle

\begin{abstract}
This paper addresses the problem of counting solutions to the ASEP-Bethe ansatz equations. In an earlier preprint, Brattain, Do, and Saenz presented an algebro-geometric method for counting the number of such solutions. As their method relies on sophisticated machinery, it requires advanced knowledge of algebraic geometry and topology.

In this paper, we provide an alternative and simpler method for counting solutions, based on a standard intersection theory on projective varieties. We take a conventional approach in algebraic geometry: ``Projectivize the problem and count degrees.''
We introduce a projectivized version of the ASEP-Bethe equations, whose solution set is realized as a subscheme of a certain projective variety. Varying the hopping-rate parameter $p$, we give a one-parameter family $\mathcal{F}$ of the solution sets.
While not all solutions are physically admissible, we can extract the essential part $\mathcal{F}_{\mathrm{ES}}\subset \mathcal{F}$ by omitting unphysical solutions. We show that, for the ASEP with $n$ sites and $m$ particles, the essential part $\mathcal{F}_{\mathrm{ES}}|_p$ for generic $p$ consists of $\binom{n}{m}$ points with multiplicity, which equals the expected number of admissible Bethe roots to the ASEP-Bethe ansatz equation.

\bigskip

Keywords: Bethe ansatz equations, Asymmetric Simple Exclusion Process, Heisenberg model, intersection theory, counting problem

MSC: 14N10, 82C22, 60K35, 60J27

\end{abstract}

\section{Introduction}

Let $m$ and $n$ be positive integers with $m<n$.
For $1\leq i\leq n$, let 
$$S_i(z_1,\dots,z_m)$$ be a symmetric polynomial in $z_1,\dots,z_m$.
In this paper, we consider the system of equations
\begin{equation}\label{eq:Bt-equation}
S_0(z_1,\dots,z_m)z_j^n-S_1(z_1,\dots,z_m)z_j^{n-1}+\dots+(-1)^n
S_n(z_1,\dots,z_m)=0,\qquad (j=1,2,\dots,m).
\end{equation}

The system \eqref{eq:Bt-equation} includes a wide class of Bethe-type equations for quantum models such as the TASEP with periodic boundary (see below), the TASEP with open boundary (Example \ref{eq:q_open_TASEP}), the periodic XXX chain (Example~\ref{eq:q_periodic_XXZ}), and the periodic ASEP (Example~\ref{ex:priodic_ASEP}).

For example, if we substitute
\[
S_i:=
\begin{cases}
g\binom{n}{i}\cdot \prod_{i=1}^mz_i & (i\neq n-m),\\
g\binom{n}{n-m}\cdot \prod_{i=1}^mz_i+1 & (i=n-m)
\end{cases}
\]
to \eqref{eq:B_P_TASEP}, we obtain the Bethe equations for the quasi-periodic TASEP~\cite{iwao2025bethe}:
\begin{equation}\label{eq:B_P_TASEP}
\frac{z_j^m}{(1-z_j)^n}=(-1)^{m+1}g\cdot \prod_{i=1}^mz_i,\qquad (j=1,2,\dots,m),
\end{equation}
where $g\neq 0$ is a parameter related to the fugacity.
If $g=1$, this system corresponds to the usual periodic TASEP.

We refer to the system of equations \eqref{eq:Bt-equation} as the \textit{ASEP-Bethe-type system}.
The following equations in Examples \ref{eq:q_open_TASEP}--\ref{ex:priodic_ASEP} are also specializations of the ASEP-Bethe-type system.

\begin{example}\label{eq:q_open_TASEP}
The Bethe equations for the open-boundary TASEP~\cite{crampe2018equivalent} are given by
\[
(\alpha z_j+p)(z_j-1)^{2n}(z_j^2-1)
=z_j^{n+1}
(\alpha z_j+1)(z_j+\beta)\prod_{k=1}^n\left(z_j-\frac{1}{z_k}\right),\quad (j=1,2,\dots,n),
\]
where $p,\alpha,\beta$ are real parameters.
\end{example}

\begin{example}\label{eq:q_periodic_XXZ}
The Bethe equations for the periodic Heisenberg XXZ spin chain~\cite{alba2013bethe,noh2000origin} is given by
\[
z_j^n=(-1)^{m-1}\prod_{k:k\neq j}\frac{z_kz_j-pz_j+1}{z_kz_j-pz_k+1},\qquad (j=1,2,\dots,m),
\]
where $p$ is a real parameter.
\end{example}

\begin{example}\label{ex:priodic_ASEP}
The Bethe equations for the quasi-periodic ASEP~\cite{brattain2015completeness,prolhac2021riemann} are given by
\begin{equation}\label{eq:ASEP_Bethe}
z_j^n=(-1)^{m-1}g\prod_{k:k\neq j}\frac{qz_jz_k-z_j+p}{qz_jz_k-z_k+p}\qquad
(j=1,2,\dots,m)
\end{equation}
where $g$ and $p$ are real parameters and $q=1-p$.
By substituting $p=1$ ($q=0$), $z_j\mapsto 1-z_j$, and $g\mapsto g^{-1}$, one obtains the periodic TASEP-Bethe equations.
\end{example}

\begin{rem}
Let $s\in \RR\setminus \{0\}$.
Putting $z_i=\frac{u_i+is}{u_i-is}$ to \eqref{eq:ASEP_Bethe} yields
\begin{equation}
\left(
\frac{u_j+is}{u_j-is}
\right)^n=(-1)^{m-1}g\prod_{k:k\neq j}\frac{qu_j-pu_k+is}{qu_k-pu_j+is}\qquad
(j=1,2,\dots,m).
\end{equation}
In particular, when $p=q=\frac{1}{2}$, $s=\frac{1}{2}$, and $g=1$, we obtain the system
\begin{equation}\label{eq:trans_Bethe}
\left(
\frac{u_j+\frac{i}{2}}{u_j-\frac{i}{2}}
\right)^n=(-1)^{m-1}\prod_{k:k\neq j}\frac{u_j-u_k+i}{u_k-u_j+i}\qquad
(j=1,2,\dots,m).
\end{equation}
which is a well-known form of the Bethe equations for the periodic XXX chain.
\end{rem}

The symmetry of the equations implies that the ordering of $z_1,z_2,\dots,z_m$ is irrelevant.
Then, in what follows, a \textit{solution} to a ASEP-Bethe type system means the (unordered) multiset $\{z_1,\dots,z_m\}$ satisfying \eqref{eq:Bt-equation}.
Such a solution is also called a \textit{Bethe solution} or a \textit{Bethe root}.

\subsection{The aim of this paper}

Counting the number of solutions to the Bethe equations, or determining them explicitly, is a long-standing problem that has been the subject of extensive research.
Complicating matters further, in the context of quantum physics, not all solutions to the Bethe-type systems are regarded as relevant.
For example, the Bethe equations for the XXX spin chain (Example \ref{eq:q_periodic_XXZ}) have certain ``inadmissible'' solutions, which are known as \textit{unphysical solutions}, that do not give rise to a \textit{Bethe vector} (that is, an eigenvector of the Markov matrix).
The issue of the presence of such solutions has been the subject of numerous previous studies~\cite{avdeev1986exceptional,kirillov2014singular,KirillovSakamoto2015,nepomechie2013algebraic,nepomechie2014twisting} and intensive research.


A solution $\{z_1,\dots,z_m\}$ to the ASEP-Bethe-type system \eqref{eq:Bt-equation} is called \textit{simple} if $z_1,\dots,z_m$ are pairwise distinct (Definition \ref{defi:quasi_simple}).
Simple solutions are usually expected to enjoy all desirable properties.
In fact, for the periodic ASEP (Example \ref{ex:priodic_ASEP}), there exists an algebraic formula for deriving a Bethe vector from a simple Bethe root
to (Eq.~\eqref{eq:formula_eigen_vec}. See also \cite[Eqs.~(1.6) and (1.7)]{brattain2015completeness}).
It is known that forcibly applying this formula to non-simple solutions yield only a zero vector. 
In this sense, the simpleness is essential for the construction.


In general, however, simple Bethe roots may be inadmissible.
A typical example is a \textit{singular solution} to the XXX-Bethe equations \eqref{eq:trans_Bethe}, which is of the form
\begin{equation}\label{eq:singular_solution}
\{u_1,\dots,u_m\}=\{\tfrac{i}{2},-\tfrac{i}{2},u_3,\dots,u_m\}.
\end{equation}
There exists an algebraic formula for the Bethe vector 
\[
B(u_1)B(u_2)\dots B(u_m)|0\rangle,
\]
where $B(u)$ is a certain operator acting on the ground state vector $|0\rangle$. (For detailed definitions of this expression, see the review article \cite[Section 4]{levkovich2016bethe}).
However, this does not work when $\{u_1,\dots,u_m\}$ is singular.
Via perturbation and normalization~\cite{nepomechie2013algebraic,nepomechie2014twisting}, some singular solutions yield a Bethe vector, while others do not.
The former are referred to as \textit{physical singular solutions}, and the latter are known as \textit{unphysical singular solutions}.

Through the variable change $z=\frac{u+\frac{i}{2}}{u-\frac{i}{2}}$, the expression \eqref{eq:singular_solution} can be transformed into
\[
\{z_1,\dots,z_m\}=\{\infty,0,z_3,\dots,z_m\}.
\]
It would be natural to say that we need some ``projectivization'' to deal with such solutions in a consistent fashion.
This paper aims to present a new algebro-geometric method for projectivizing the problem and distinguishing which solutions are relevant and which are irrelevant.


\subsection{Organization of this paper}

This paper is organized as follows.
In Section \ref{sec:project}, we introduce the projectivized Bethe system \eqref{eq:induced_system}, which is a system of equations defined over the projective variety $\PP^m\times \PP^n$.
This system captures all \textit{quasi-simple solutions} (Definition \ref{defi:quasi_simple}) to the ASEP-Bethe-type system \eqref{eq:Bt-equation}.
Intuitively, quasi-simple solutions are ``perturbation-friendly'' solutions, which are ``very near to'' simple solutions (see \eqref{eq:zero_system_ep}).
We expect quasi-simple solutions to exhibit desirable properties.
Our framework has the advantage of automatically filtering out ``non-perturbable'' solutions.

In Section \ref{sec:app_to_ASEP}, we introduce the projectivizetized ASEP-Bethe system \eqref{eq:system_a_and_b_m}, which captures all quasi-simple solutions to the ASEP-Bethe equations~(Example \ref{ex:priodic_ASEP}). 
Varying the hopping-rate-parameter $p\in \CC$ yields a one-parameter family $S$ of quasi-simple solutions (\S \ref{sec:def_of_S}).
As $S$ is a generically $1$-dimensional scheme, we can unambiguously determine the direction of the perturbation along it if $S$ is non-singular.
(If $S$ is singular at a specified point $P$, one need to work on the desingularization $t_P:\widetilde{S}\to S$.
In such a case, there are finitely many directions of the perturbation according to the branches $Q_1,\dots,Q_d\in t_P^{-1}(P)$.)

Through direct analysis (see \S \ref{sec:Irr_of_S} and \S \ref{sec:special_fiber}) on $S$, we show that
\begin{enumerate}
\item[(i)] irreducible components of $S$ are classified into four types (type X, $\infty$, E, and $Z_{\mathrm{steady}}$) according to the way to intersect with $\{p=0\}$ (see Figure \ref{fig:1}), and
\item[(ii)] the total number of points in all irreducible components of type E and $Z_{\mathrm{steady}}$ is $\binom{n}{m}$ with multiplicity, which coincides with the expected number of admissible Bethe roots (that is, solutions to the ASEP-Bethe equations that yield a non-zero Bethe vector).
\end{enumerate}
We exclude the components of types X and $\infty$ from our argument\footnote{Precisely speaking, if an irreducible component of type X or $\infty$ intersects other irreducible component of type E or $Z_{\mathrm{steady}}$, an admissible solution can be obtained from the intersection.
Such points cause no difficulty in the subsequent discussion, as they can be reconstructed from type E or $Z_{\mathrm{steady}}$ component by taking a continuous limit.
} because the explicit formula \eqref{eq:formula_eigen_vec} for the Bethe vector is not applicable for generic points in type $\infty$ component, and type $X$ component exists only for finitely many $p$.

Our results provide strong evidence for the \textit{completeness}, which claims that the number of admissible solutions to the ASEP-Bethe equations is $\binom{n}{m}$. 
On the other hand, for a complete proof of this statement, we must show further that ``all points in the component of types E and $Z_{\mathrm{steady}}$ yield a non-zero Bethe vector'' and ``these Bethe vectors are linearly independent.''

Because our method does not directly make any claims about Bethe vectors, it does \textit{not} provide a complete proof of the completeness.
However, at least we can say that:
\begin{itemize}
\item any limit of simple solutions is quasi-simple (= the method of perturbation is not applicable for non-quasi-simple solutions) (Section \ref{sec:project}), 
\item for generic $p$, the number of quasi-simple solutions to the ASEP-Bethe equations, contained in the finite part $\CC\subsetneq \PP^1$, is at most $\binom{n}{m}$, and
\item for generic $g$ and $p$, there exists exactly $\binom{n}{m}$ quasi-simple solutions to the ASEP-Bethe equations contained in $\CC\subsetneq \PP^1$.
\end{itemize}

In Section \ref{sec:future}, we consider the reduction to the special case $p=q=\frac{1}{2}$ and $g=1$ to demonstrate the effectiveness of our method.
In this case, the ASEP reduces to the \textit{SEP} (Symmetric Exclusion Process).
Under the reduction to the SEP, some functions such as a component in Bethe vectors may diverge by simple substitution $p=q=\frac{1}{2}$.
In such cases, perturbing the parameters is effective~\cite{nepomechie2013algebraic,nepomechie2014twisting}. 
While existing methods required a careful choice of the perturbation direction, in our approach it is determined automatically because $S$ is $1$-dimensional.
At the end of this paper, we present a method for actually obtaining solutions via perturbation for a concrete example with $(m,n)=(2,4)$.

\subsection{Related works}

This work is inspired by two earlier works: One is the preprint~\cite{brattain2015completeness} by Brattain-Do-Saenz, which introduced a rigorous algebro-geometric method for counting Bethe roots for the periodic ASEP. 
The other is \cite{jiang2018algebraic} by Jiang-Zhang, which employs a computational method based on symmetric polynomials to calculate Bethe roots.
This paper aims to present a ``easy-to-access'' formulation of the counting problem by reformulating Jiang-Zhang's symmetric polynomial approach within an algebro-geometric framework, and combining it with intersection theory.

For the TASEP (Totally Asymmetric Simple Exclusion Process), which is obtained from the ASEP by substituting $p=0$, the counting problem has been extensively studied in the literature. 
Of particular relevance to the present study are Prolhac's serial works~\cite{prolhac2020riemann,prolhac2020riemannKPZ,prolhac2021riemann} on a Riemann surface approach.
Building on Prolhac's method, Motegi and the author of the present paper provided a proof that the number of solutions to the TASEP-Bethe equations is $\binom{n}{m}$ with multiplicity~\cite{iwao2025bethe}.
They also show that, for a generic $g$, the number of solutions is exactly $\binom{n}{m}$.
Although the TASEP is a special case of the ASEP, this fact does not imply that the analysis on TASEP is simpler than that on the ASEP. 
Their method using Riemann surfaces reflects deep structures, such as the symmetry of the system, and warrants further investigation.

\section{The ASEP-Bethe-type equation and its projectivization}\label{sec:project}

\subsection{Projectivized Bethe system}\label{sec:alg_form}

Let $z=\{z_1,\dots,z_m\}$ be a (unordered) set of indeterminates and
\[
e_k=e_k(z)=\sum_{1\leq i_1<\dots<i_k\leq m} \prod_{j=1}^kz_{i_j} 
\]
be the $k$-th elementary symmetric polynomial.
We set $e_0=1$ and $e_k=0$ for $k>m$.
Then, for each $j$, we have \textit{Vieta's formula}: 
\begin{equation}
z_j^m-e_1z_j^{m-1}+\dots+(-1)^me_m=0.
\end{equation}

Consider the ASEP-Bethe-type system \eqref{eq:Bt-equation}.
Without loss of generality, we may assume the coefficient of the highest term in \eqref{eq:Bt-equation} to be $1$, by subtracting the trivial equation $z_j^{n+1}-e_1z_j^{n}+\dots+(-1)^me_mz^{n-m+1}=0$ from it if necessary.
Hence, we shall consider the system
\begin{equation}\label{eq:monic_Bt-equation}
z_j^n-S_1(z_1,\dots,z_m)z_j^{n-1}+\dots+(-1)^n
S_n(z_1,\dots,z_m)=0,\qquad (j=1,2,\dots,m).
\end{equation}

Since each $S_i(z) = S_i(z_1, \dots, z_m)$ is a symmetric polynomial, fixing a (unordered) set of $m$ points $\{z_1, \dots, z_m\}\subset \CC$ determines the values $S_i$ accordingly. 
Then, the system of equations \eqref{eq:monic_Bt-equation} is equivalent to the following inclusion relationship of sets:
\begin{equation}\label{eq:incl_sets}
\left\{
z_1,\dots,z_m
\right\}\subset 
\left\{
\lambda\in \CC:
\begin{aligned}
\lambda^n - S_1(z)\lambda^{n-1} + \cdots + (-1)^nS_n(z) = 0
\end{aligned}
\right\}.
\end{equation}

This interpretation can be projectivized simply by replacing the complex line $\CC$ with the projective line $\PP^1$. 
Let $(\PP^1)^{[n]} = \mathrm{Hilb}^n(\PP^1)$ be the Hilbert scheme \cite[\S 7.1]{fantechi2005fundamental} of $n$ points on $\PP^1$. 
We identify a (unordered) set $\{z_1, \dots, z_n\}\subset \PP^1$ (with multiplicities) with the polynomial
\[
a_0z^n - a_1z^{n-1} + \cdots + (-1)^na_n,
\]
whose roots are $\{z_1,\dots,z_n\}$.
In this way, an $n$-point configuration with multiplicity can be associated with an element $[a_0 : a_1 : \cdots : a_n]\in \PP^n$. 
This correspondence gives an isomorphism of schemes $(\PP^1)^{[n]} \simeq \PP^n$ (see \cite[Example 7.1.3]{fantechi2005fundamental}).

Under this identification, the $\cup$-operator
\[
\cup:(\{z_1,\dots,z_m\},\{w_1,\dots,w_{n-m}\})\mapsto \{z_1,\dots,z_m,w_1,\dots,w_{n-m}\}
\]
is realized by the map
\[
\begin{gathered}
Q:\PP^{m}\times \PP^{n-m}\to \PP^n;\\
\left([a_0:\dots:a_m],[b_0:\dots:b_{n-m}]\right)\mapsto 
\left(
a_0b_0:\sum_{i+j=1}a_ib_j :
\sum_{i+j=2}a_ib_j:\dots:a_mb_{n-m}
\right).
\end{gathered}
\]

Define the injection $\iota$ by
\[
\iota:\PP^m\times \PP^{n-m}\to \PP^m\times \PP^{n};\quad
(A,B)\mapsto \left(A,Q(A,B)\right).
\]
Then, the image
\[
Z_{m,n}:=\iota(\PP^m\times \PP^{n-m})
\]
parameterizes all nested sets
\begin{equation}\label{eq:nested_set}
Z\subset \Lambda \subset \PP^1
\end{equation}
satisfying $|Z|=m$ and $|\Lambda|=n$.
Let
$(e,E)=([e_0:e_1:\dots:e_m],[E_0:E_1:\dots:E_n])$
be the homogeneous coordinate of $\PP^m\times \PP^n$.
Then, for $(e,E)\in Z_{m,n}$, the nested set \eqref{eq:nested_set} is characterized by
\[
\begin{aligned}
\Lambda&=\{[\lambda_0:\lambda_1]\in \PP^1:E_0\lambda_1^n-E_1\lambda_1^{n-1}\lambda_0+\dots+(-1)^nE_n\lambda_0^n=0 \},\\
Z&=\{[\lambda_0:\lambda_1]\in \PP^1:e_0\lambda_1^m-e_1\lambda_1^{m-1}\lambda_0+\dots+(-1)^me_m\lambda_0^m=0 \}.
\end{aligned}
\]

Let $\widehat{S}_k(e_1, \dots, e_n)$ denote the polynomial obtained by expressing $S_k(z_1, \dots, z_m)$ in terms of the elementary symmetric polynomials $e_1, \dots, e_n$. We write its homogenization as 
\[
\mathcal{S}_k(e_0, e_1, \dots, e_m) := e_0^{\max[\deg(\hat{S}_k),1]} \cdot \hat{S}_k(e_1/e_0, \dots, e_n/e_0).
\]
Under the above setting, the inclusion relation \eqref{eq:incl_sets} can be characterized by the following system of equations defined on $Z_{m,n}$:
\begin{equation}\label{eq:system_on_Z}
\frac{E_k}{E_0}=\frac{\mathcal{S}_k(e_0,e_1,\dots,e_m)}{(e_0)^{
\max[\deg(\hat{S}_k),1]}},\quad (k=1,2,\dots,n).
\end{equation}

Let $\mathrm{Sol}$ be the solution set
\begin{equation}\label{eq:set_incl_affine}
\mathrm{Sol}:=
\left\{
([e_0:\dots:e_m],[E_0:\dots:E_n]):
\frac{E_k}{E_0}=\frac{\mathcal{S}_k(e_0,e_1,\dots,e_m)}{(e_0)^{
\max[\deg(\hat{S}_k),1]}},\
(k=1,2,\dots,n)
\right\}\cap Z_{m,n}.
\end{equation}
Once the homogeneous coordinates $(e,E)\in \mathrm{Sol}$ are obtained, a solution to the ASEP-Bethe-type equation \eqref{eq:Bt-equation} can be obtained by solving the algebraic equation
\[
e_0z^m-e_1z^{m-1}+\dots+(-1)^me_m=0.
\]

Through the pullback
\[
\iota^\ast:e_k\mapsto a_k,\quad E_k\mapsto \sum_{i+j=k}a_ib_j,
\]
of the birational map $\iota$, we induce from \eqref{eq:system_on_Z} the system of equations
\begin{equation}\label{eq:induced_system}
a_0^{
[\deg(\hat{S}_k)-1]_+
}\cdot
\sum_{i+j=k}a_ib_j= b_0\cdot \mathcal{S}_k(a_0,a_1,\dots,a_m),\qquad (k=1,2,\dots,n).
\end{equation}
Here, $[x]_+:=\max[x,0]$ and $(a,b):=([a_0:a_1:\dots:a_m],[b_0:b_1:\dots:b_{n-m}])$ are the homogeneous coordinates of $\PP^m\times \PP^{n-m}$.

We will call \eqref{eq:induced_system} as the \textit{projectivized Bethe system}.
From a solution to \eqref{eq:induced_system}, we can derive a solution to the ASEP-Bethe-type system~\eqref{eq:monic_Bt-equation} by solving the algebraic equation
\[
a_0z^m-a_1z^{m-1}+\dots+(-1)^ma_m=0.
\]

If the leading term $S_0$ in \eqref{eq:Bt-equation} is not $1$, we subtract the trivial equation $z_j^{n+1}-e_1z_j^{n}+\dots+(-1)^{m}e_mz_{n-m+1}=0$ from \eqref{eq:Bt-equation} to obtain the monic system
\[
z_j^{n+1}-e_1z_j^{n}+\dots+(-1)^{m}e_mz_{n-m+1}
-\sum_{k=0}^n(-1)^kS_k(z_1,\dots,z_m)z_j^{n-k}=0.
\]
Its projectivized system \eqref{eq:induced_system} is expressed as
\[
a_0^{[\mathrm{deg}(\hat{S}_k)-1]_+}\cdot
\sum_{i+j=k}a_ib_j= b_0\cdot \left\{a_0^{[\mathrm{deg}(\hat{S}_k)-1]_+}\widetilde{a}_k+
\mathcal{S}_k(a_0,a_1,\dots,a_m)\right\},\qquad (k=1,2,\dots,n+1),
\]
where $\widetilde{a}_k=\begin{cases}
a_k& (0\leq k\leq m),\\
0 & (k>m).
\end{cases}$
By simplifying it, we obtain the system of equations
\begin{equation}\label{eq:modif_proj}
a_0^{[\mathrm{deg}(\hat{S}_{k-1})-1]_+}\cdot
\sum_{\substack{i+j=k\\i\geq 0,\ j\geq 1}}a_ib_j= b_0\cdot \mathcal{S}_{k-1}(a_0,a_1,\dots,a_m),\qquad (k=1,2,\dots,n+1),
\end{equation}
which is also refereed to as the \textit{projectivized Bethe system}. 

The affine part of the solution set of \eqref{eq:induced_system} is in a one-to-one correspondence with quasi-simple solutions to the original system \eqref{eq:Bt-equation}.
However, the solution set of \eqref{eq:induced_system} itself need not coincide with any compactification of that of \eqref{eq:Bt-equation}.
For example, when $(m,n)=(1,2)$ and $S_1(z)=S_2(z)=z^2$, the ASEP-Bethe-type system
\begin{equation}\label{eq:original_1}
z^2-(z^2)\cdot z^1+(z^2)=0
\end{equation}
admits solutions $z=0,2$, while the projectivized system \eqref{eq:induced_system} admits
\[
([a_0:a_1],[b_0:b_1])=([1:0],[1:0]),([1:2],[1:2]),([0:1],[0:1]).
\]
The first and second solutions correspond to $z=0$ and $z=2$, respectively.
However, the third solution does not corresponds to a solution nor a limit of solutions to \eqref{eq:original_1}.
This example demonstrates that our projectivization is distinct from a mere compactification of the solution set.

\subsection{Quasi-simple solutions and perturbation}

When dealing with inclusion relations such as \eqref{eq:incl_sets}, we must be careful about the difference between \textit{sets} and \textit{multisets}.
Indeed, while the relation \eqref{eq:incl_sets} represents an inclusion of \textit{sets}, the variety $Z_{m,n}$ parametrizes the inclusion of \textit{multisets}. 
Strictly speaking, a solution $z=\{z_1,\dots,z_m\}$ to \eqref{eq:incl_sets} is in a one-to-one correspondence to a point in \eqref{eq:set_incl_affine} \textit{only if $z_1, \dots, z_m$ are all distinct}. 
To clarify these differences, we introduce the following definitions.

\begin{defi}[Simple and quasi-simple solutions]\label{defi:quasi_simple}
Let $z=\{z_1,\dots,z_m\}\subset \CC$ be a solution to the ASEP-Bethe-type system \eqref{eq:monic_Bt-equation}. 
(Hence, $z$ satisfies the inclusion relation \eqref{eq:incl_sets} of sets.)
Then, we will say that:
\begin{enumerate}
\item[(1)] $z$ is a simple solution if $z_1,\dots,z_m$ are all distinct, and
\item[(2)] $z$ is a quasi-simple solution if it satisfies the inclusion relation \eqref{eq:incl_sets} as multisets.
\end{enumerate}
Obviously, every simple solution is quasi-simple.
\end{defi}

We can say that $\mathrm{Sol}$ (Eq.~\eqref{eq:set_incl_affine}) is a parameter space of all quasi-simple solutions to the ASEP-Bethe-type system.

Before going further, we give three examples illustrating the differences between simple, quasi-simple, and non-quasi-simple solutions.

\begin{example}[A quasi-simple solution that is not simple]
Consider the ASEP-Bethe-type system
\begin{equation}\label{eq:zero_system}
\begin{cases}
z_1^3-\frac{3}{2}(z_1+z_2)z_1+2=0,\\
z_2^3-\frac{3}{2}(z_1+z_2)z_2+2=0.
\end{cases}
\end{equation}
Here, the coefficients are
\[
(\widehat{S}_0,\widehat{S}_1,\widehat{S}_2,\widehat{S}_3)=(1,0,-\frac{3}{2}e_1,-2).
\]
Then, the projectivized problem (Eq.~\eqref{eq:induced_system}) is
\[
\begin{cases}
a_0b_1+a_1b_0=0,\\
a_1b_1+a_2b_0=-\frac{3}{2}a_1b_0,\\
a_2b_1=-2a_0b_0,\\
\end{cases}\qquad
([a_0:a_1:a_2],[b_0:b_1])\in \PP^2\times \PP^1.
\]
It is straightforward to obtain the three solutions
\[
([a_0:a_1:a_2],[b_0:b_1])=([1:2:1],[1:-2]),\quad
\left(
\left[1:-\zeta_{\pm}:-2\zeta_{\pm}^{-1}\right],\left[
1:\zeta_{\pm}\right]
\right),
\]
where $\zeta_+$ and $\zeta_{-}$ are two roots of $2\zeta^2-\zeta+2$.
By solving $a_0z^2-a_1z+a_2=0$, we obtain the three solutions:
\[
\{z_1,z_2\}=\{1,1\},\quad \left\{
\eta_{+1},\eta_{+2}
\right\},\quad \{\eta_{-1},\eta_{-2}\},
\]
where $\eta_{\pm 1}$ and $\eta_{\pm 2}$ are two roots of $\eta^2+\zeta_\pm\eta -2\zeta_\pm^{-1}$.
The first solution $\{1,1\}$ is not simple, but quasi-simple, because the following multiset inclusion holds:
\[
\{1,1\}\subset \left\{\lambda\in \CC:
a_0\lambda^3-\frac{3}{2}a_1\lambda+2a_0=
\lambda^3-3\lambda+2=0\right\}=\{1,1,-2\}.
\]
\end{example}

\begin{example}[Non-quasi-simple solutions]
Assuming $z_1 = z_2$ in the system \eqref{eq:zero_system}, we also obtain two additional solutions:
\[
\{z_1,z_2\}=\{1+\sqrt{3},1+\sqrt{3}\},\quad \{1-\sqrt{3},1-\sqrt{3}\}.
\]
However, these are not quasi-simple.
In fact, for example, the first solution corresponds to the inclusion relation
\[\begin{aligned}
\{1+\sqrt{3},1+\sqrt{3}\}
&\subset 
\left\{
\lambda\in \CC: \lambda^3-3(1+\sqrt{3})\lambda+2=0
\right\}\\
&=\left\{
1+\sqrt{3},
\frac{-1-\sqrt{3}+ \sqrt{6\sqrt{3}}}{2},
\frac{-1-\sqrt{3}- \sqrt{6\sqrt{3}}}{2}
\right\},
\end{aligned}
\]
which is not a multiset inclusion.
\end{example}


\begin{example}[Non-monic system]
Let $m=2$ and $n=3$.
Consider the ASEP-Bethe-type system
\begin{equation}\label{eq:second_system}
\begin{cases}
2(z_1z_2)z_1^3-z_1+3(z_1+z_2)=0,\\
2(z_1z_2)z_2^3-z_2+3(z_1+z_2)=0.
\end{cases}
\end{equation}
The coefficients $\widehat{S}_k$ are given by
\[
(\widehat{S}_0,\widehat{S}_1,\widehat{S}_2,\widehat{S}_3)
=(2e_2,0,-1,-3e_1).
\]
Since $\widehat{S}_0$ is not $1$, we use \eqref{eq:modif_proj} to derive the projectivized system
\[
\begin{cases}
a_0b_1=2a_2b_0,\\
a_1b_1+a_0b_2=0,\\
a_2b_1+a_1b_2=-a_0b_0,\\
a_2b_2=-3a_1b_0.
\end{cases}
\]
Direct calculations derive the six solutions
\[
([a_0:a_1:a_2],[b_0:b_1:b_2])=
([\pm\sqrt{-2}:0:1],[1:\mp\sqrt{-2}:0]),\quad
([\zeta^2:-2\zeta:3],[1:\zeta^2:2\zeta]),
\]
where $\zeta$ is an arbitrary complex number satisfying $\zeta^4=6$.
The corresponding solutions $\{z_1,z_2\}$ are given by
\[
\{z_1,z_2\}=
\{2^{-\frac{1}{4}}e^{\frac{\pi}{4}i},
2^{-\frac{1}{4}}e^{\frac{5\pi}{4}i}\},\ 
\{2^{-\frac{1}{4}}e^{\frac{3\pi}{4}i},
2^{-\frac{1}{4}}e^{\frac{7\pi}{4}i}\},\ 
\{(-1+\sqrt{-2})\zeta^{-1},(-1-\sqrt{-2})\zeta^{-1}\}.
\]
In addition, \eqref{eq:second_system} admits the non-quasi-simple solutions
\[
\{z_1,z_2\}=\{0,0\},\quad \{\eta,\eta\},\quad \text{where}\quad \eta^4=-\frac{5}{2}.
\]
\end{example}

Intuitively, quasi-simple solutions are ``perturbation-friendly.''
As a motivating example, consider the perturbed system 
\begin{equation}\label{eq:zero_system_ep}
\begin{cases}
z_1^3-\frac{3}{2}(z_1+z_2)z_1+2=\epsilon,\\
z_2^3-\frac{3}{2}(z_1+z_2)z_2+2=\epsilon
\end{cases}
\end{equation}
of \eqref{eq:zero_system_ep}, where $\epsilon$ is a small parameter.
If $|\epsilon|$ is sufficiently small,
the quasi-simple solution $\{z_1,z_2\}=\{1,1\}$ is perturbed as 
\[
\{z_1,z_2\}=\left\{1+\frac{\epsilon^{1/2}}{2}-\frac{\epsilon}{12}+\cdots,1-\frac{\epsilon^{1/2}}{2}-\frac{\epsilon}{12}+\cdots\right\},
\]
which is simple.
On the other hand, however, the non-quasi-simple solution $\{z_1,z_2\}=\{1+\sqrt{3},1+\sqrt{3}\}$ is perturbed as
\[
\{z_1,z_2\}=\left\{
1+\sqrt{3}+\frac{1}{6}\epsilon+\cdots,
1+\sqrt{3}+\frac{1}{6}\epsilon+\cdots
\right\},
\]
which is not simple.
In this sense, the first solution can be said to be ``very close to'' simple solutions, while the second is not.

In fact, all limits of simple solutions must be quasi-simple.
This fact follows from the continuity of roots for the nested sets \eqref{eq:incl_sets}.
%
This is why we should disregard non-quasi-simple solutions and focus on quasi-simple ones.

\begin{rem}
This statement does not guarantee that ``all quasi-simple solutions can be made simple by a given perturbed system.''
Whether simple solutions can be obtained in a similar manner for an individual system requires a case-by-case analysis.
\end{rem}


If the projectivized Bethe system admits a solution satisfying
\[
a_0=a_1=\dots=a_{k-1}=0,\quad a_k\neq 0,\qquad \text{for some $k\geq 1$,}
\]
it corresponds to a solution to the ASEP-Bethe-type system that contains $\infty\in \PP^1$:
\[
\{z_1,\dots,z_{m-k},\overbrace{\infty,\infty,\dots,\infty}^{k}\},
\]
where $z_1,\dots,z_{m-k}$ are the roots of $a_{k}z^{m-k}-a_{k+1}z^{m-k-1}+\dots+(-1)^{m-k}a_m$.

\begin{example}[Bethe roots at infinity]
Let $m=3$ and $n=4$.
Consider the ASEP-Bethe-type system
\begin{equation}
\begin{cases}
(z_1+z_2+z_3)z_1^4+(z_1z_2z_3)z_1^2-(z_1z_2+z_1z_3+z_2z_3)z_1=0,\\
(z_1+z_2+z_3)z_2^4+(z_1z_2z_3)z_2^2-(z_1z_2+z_1z_3+z_2z_3)z_2=0,\\
(z_1+z_2+z_3)z_3^4+(z_1z_2z_3)z_3^2-(z_1z_2+z_1z_3+z_2z_3)z_3=0,
\end{cases}
\end{equation}
whose projectivization is
\[
\begin{cases}
a_0b_1=a_1b_0,\\
a_1b_1+a_0b_2=0,\\
a_2b_1+a_1b_2=a_3b_0,\\
a_3b_1+a_2b_2=a_2b_0,\\
a_3b_2=0.
\end{cases}
\]
Here, $(a,b)=([a_0:a_1:a_2:a_3],[b_0:b_1:b_2])$ is the homogeneous coordinate of $\PP^3\times \PP^2$.
This projectivized system has six solutions:
\[
\begin{aligned}
(a,b)
=
&([1:0:0:0],[1:0:0]),\quad
([1:\pm i:-1:0],[1:\pm i:1]),\\
&([0:0:1:0],[1:0:1]),\quad
([0:0:1:\pm 1],[1:\pm 1:0]),
\end{aligned}
\]
which correspond to the Bethe roots
\[
\{z_1,z_2,z_3\}=
\{0,0,0\},\quad 
\{0,\frac{\sqrt{3}\pm \sqrt{-1}}{2},\frac{-\sqrt{3}\pm \sqrt{-1}}{2}\},\quad
\{0,\infty,\infty\},\quad
\{\pm 1,\infty,\infty\}.
\]
\end{example}

\section{Application to ASEP}\label{sec:app_to_ASEP}

\subsection{Preliminaries: Bethe vectors}

The ASEP (Asymmetric Simple Exclusion Process) is a Markov process in which particles hop randomly between neighboring sites of a one-dimensional lattice  (Figure \ref{fig:ASEP}). 
Let $n$ denote the size of the system, and $m$ denote the number of particles.
Each site is either occupied by a particle or empty. In an infinitesimal time interval $[t, t+dt]$, each particle hops to the neighboring site on the right with probability $p\,dt$ and to the one on the left with probability $q\,dt$, provided that the target site is empty.

\begin{figure}[htbp]
\centering
\begin{tikzpicture}
\def\haba{.8} 
\def\kosuu{14}
\draw[thick](0,0)--({\kosuu * \haba},0);
\draw[thick](0,\haba)--({\kosuu * \haba},\haba);
\foreach \i in {0,1,2,...,\kosuu}
{
\draw[thick]({\i * \haba},0)--++(0,\haba) ;
}
\foreach \i in {1,2,5,7}
{
\draw (\haba*\i, 0)++({\haba/2}, {\haba/2}) circle[radius = 0.4*\haba];
\fill[gray, fill opacity=.5] (\haba*\i, 0)++({\haba/2}, {\haba/2}) circle[radius = 0.4*\haba];
}
\draw[->, >=stealth, thick] (\haba*5, 0)++({\haba/2+.1}, {\haba}) ..controls ++({\haba/2}, {\haba})..node[above]{$p$}
++({\haba}, 0);
\draw[->, >=stealth, thick] (\haba*5, 0)++({\haba/2-.1}, {\haba}) ..controls ++({-\haba/2}, {\haba})..node[above]{$q$}
++({-\haba}, 0);
\draw[->, >=stealth, thick] (\haba*2, 0)++({\haba/2+.1}, {\haba}) ..controls ++({\haba/2}, {\haba})..node[above]{$p$}
++({\haba}, 0);
\draw[->, >=stealth,thick] (\haba*2, 0)++({\haba/2-.1}, {\haba}) ..controls ++({-\haba/2}, {\haba})..
++({-\haba}, 0);
\draw[thick](\haba*2, 0)++({-\haba/2+.2}, {\haba+.3})++(-.1,0)--++(.4,.4);
\draw[thick](\haba*2, 0)++({-\haba/2+.6}, {\haba+.3})++(-.1,0)--++(-.4,.4);
\end{tikzpicture}
\caption{A typical configuration of the ASEP.
The number of sites is $n=14$ and the number of particles is $m=4$.
}
\label{fig:ASEP}
\end{figure}

Let 
\[
\Omega:=\left\{(x_1,\dots,x_m)\in \ZZ^m:
1\leq x_1<\dots<x_m\leq n
\right\}
\]
be the phase space, where $x=(x_1,\dots,x_n)\in \Omega$ denotes the state in which the $i$-th particle is at position $x_i$ for all $i$.
Then, the occurrence probability of the state $x\in \Omega$ at time $t$ is determined by the \textit{master equation} (see, for example, \cite[Section 1]{golinelli2006asymmetric})
\begin{equation}\label{eq:master_eq}
\begin{aligned}
\frac{\partial}{\partial t}\varphi(x;t)
=\sum_{i=1}^m 
\left\{
\begin{aligned}
&
p\varphi(x-\mathbf{e}_i;t)\cdot\mathbf{1}_{\{x_{i}-x_{i-1}>1\}}
+q\varphi(x+\mathbf{e}_i;t)\cdot\mathbf{1}_{\{x_{i+1}-x_{i}>1\}}\\
&\qquad -
p\varphi(x;t)\cdot\mathbf{1}_{\{x_{i+1}-x_{i}>1\}}-q\varphi(x;t)\cdot\mathbf{1}_{\{x_{i}-x_{i-1}>1\}}
\end{aligned}
\right\},
\end{aligned}
\end{equation}
where $\varphi(x;t)\in \CC$ is a differentiable function in $t$ and $\mathbf{e}_i=(0,\dots,\overbrace{1}^i,\dots,0)$.

Under the quasi-periodic condition
\[
\varphi_g(x_2,\dots,x_m,x_1+L;t)=g\cdot \varphi_g(x_1,\dots,x_m;t),\quad (g\neq 0),
\]
the master equation \eqref{eq:master_eq} is equivalent to the matrix differential equation
\[
\frac{\partial}{\partial t}\varphi_g(t)=M_g\varphi_g(t),
\]
where $M_g\in \CC^{\binom{n}{m}\times \binom{n}{m}}$ is the \textit{Markov matrix} or the \textit{transfer matrix}, and $\varphi_g(t)=(\varphi_g(x;t))_{x\in \Omega}$ is the right probability vector.

To solve the master equation, we need to solve the spectral problem
\begin{equation}\label{eq:spec_prob}
M_gv=\lambda v,\qquad (\lambda\in \CC,\, v\in \CC^{\binom{n}{m}}).
\end{equation}
While it is difficult to solve this eigenvalue problem directly, the eigenvalues and eigenvectors can be constructed by employing the \textit{Bethe ansatz}.

For a solution $\{z_1,\dots,z_m\}$ to the ASEP-Bethe equations \eqref{eq:ASEP_Bethe}, the Bethe vector $v(z_1,\dots,z_m):=(v(x;z_1,\dots,z_m))_{x\in \Omega}$ is defined as follows~\cite[\S 1]{brattain2015completeness}:
\begin{equation}\label{eq:formula_eigen_vec}
v(x;z)=\sum_{\sigma\in S_{m}} A_\sigma(z) \prod_{i=1}^mz_{\sigma(i)}^{x_i},\qquad A_\sigma(z)=\prod_{(i,j):i<j \text{ and } \sigma^{-1}(i)>\sigma^{-1}(j)}-\frac{qz_iz_j-z_j+p}{qz_iz_j-z_i+p}.
\end{equation}
It is known that the vector $v(z)$ satisfy the equation $M_gv(z)=\lambda_{g,z}v(z)$ with eigenvalue
\begin{equation}\label{eq:eigenvalue}
\lambda_{g,z}=\sum_{i=1}^m (qz_i-1+pz_i^{-1}).
\end{equation}

This does not imply that the spectral problem \eqref{eq:spec_prob} has been solved. 
There are two major obstacles: First, it is unclear whether the ASEP-Bethe equations \eqref{eq:ASEP_Bethe} admit a sufficient number of solutions; and second, the expression \eqref{eq:formula_eigen_vec} may fail to yield a nonzero vector.
In particular, if $z$ is not a simple solution, we have $v(z)=0$ whenever $v(z)$ is well-defined.

\subsection{Projectivized ASEP-Bethe system}\label{sec:ASEP_projec}

In the sequel, we always assume
\[
m>1
\]
without otherwise stated.
(This is because the analysis on the one-particle system requires numerous exceptions.)

We now apply our algebraic method to the ASEP-Bethe equations \eqref{eq:ASEP_Bethe}.
By multiplying the denominator, we obtain the system of algebraic equations
\begin{equation}\label{eq:ASEP_poly}
z_j^n\prod_{k:k\neq j}(qz_jz_k-z_k+p)
+
(-1)^{m}g\prod_{k:k\neq j}(qz_jz_k-z_j+p)=0\qquad (j=1,\dots,m).
\end{equation}

Set
$
F(z,w):=qzw-w+p
$.
Then, \eqref{eq:ASEP_poly} is rewritten as 
\[
z_j^n\prod_{k:k\neq j}F(z_j,z_k)
+
(-1)^{m}g\prod_{k:k\neq j}F(z_k,z_j)=0\qquad (j=1,\dots,m).
\]
Define the polynomials $\Delta_j$ and $\Delta_j^T$ by
\[
\Delta_j:=\prod_{k:k\neq j} F(z_j,z_k),\qquad 
\Delta^T_j:=
\prod_{k:k\neq j}F(z_k,z_j)
\]
and
\begin{equation}\label{eq:def_of_Delta}
\Delta(z_1,\dots,z_m) := 
\prod_{j=1}^m\Delta_j
=\prod_{j=1}^m\Delta^T_j
=\prod_{(j,k):k\neq j} F(z_j,z_k).
\end{equation}
As \eqref{eq:def_of_Delta} is the product of all denominators in \eqref{eq:ASEP_Bethe}, the ASEP-Bethe equations \eqref{eq:ASEP_Bethe} and the system of equations \eqref{eq:ASEP_poly} are equivalent provided that $\Delta(z_1,\dots,z_m) \neq 0$.
We often refer to \eqref{eq:ASEP_poly} as the \textit{polynomial form} of the ASEP-Bethe equations \eqref{eq:ASEP_Bethe}.
Any solution $\{z_1,\dots,z_m\}$ to \eqref{eq:ASEP_poly} satisfying $\Delta(z_1,\dots,z_m)\neq 0$ also satisfies \eqref{eq:ASEP_Bethe}.
\begin{defi}
A solution $\{z_1,\dots,z_m\}$ to the polynomial form \eqref{eq:ASEP_poly} is called \textit{indefinite} if it satisfies $\Delta(z_1,\dots,z_m)=0$.
\end{defi}

Indefinite solutions to \eqref{eq:ASEP_poly} may or may not be solutions to \eqref{eq:ASEP_Bethe}.
However, as explained below, only the \textit{steady state solution} $z_1=z_2=\dots=z_m=1$ is relevant to us among indefinite solutions.
(The steady state solution is the sole exception and requires special handling.)

By multiplying \eqref{eq:ASEP_poly} for all $j=1,\dots,m$, we obtain
\begin{equation}\label{eq:reg_cond}
(\prod_{i=1}^mz_i^n-g^m)\cdot \Delta(z_1,\dots,z_m)=0.
\end{equation}
Therefore, for a non-indefinite solution $\{z_1,\dots,z_m\}$, it holds that
\begin{equation}\label{eq:e_m_is_const}
\prod_{i=1}^mz_i^n=g^m.
\end{equation}

For a fixed index $j$, let
\[
e_k':=e_k(z_1,\dots,z_{j-1},z_{j+1},\dots,z_m)
\]
be the elementary symmetric polynomial in the $(m-1)$ variables $z_1,\dots,z_{j-1},z_{j+1},\dots,z_m$.
The function $e_k'$ can be expressed as a polynomial in $e_1,\dots, e_m$ and $z_j^{\pm 1}$ as
\[
\begin{aligned}
e_k'&=e_k-z_je_{k-1}+z_j^2e_{k-2}-\dots+(-1)^{k}z_j^{k}e_0,\\
e_k'&=z_j^{-1}e_{k+1}-z_j^{-2}e_{k+2}+\dots+(-1)^{m-k}z_j^{-m+k}e_m.
\end{aligned}
\]
By using these expressions, $\Delta_j$ and $\Delta_j^T$ can be also expressed as a polynomial in $z_j$ with coefficients in $\ZZ[e_1,\dots,e_m,p]$.
In fact, we have
\[
\begin{aligned}
\Delta_j
&=
\prod_{k:k\neq j}\{(qz_j-1)z_k+p\}\\
&=
\sum_{i=0}^{m-1} (qz_j-1)^ip^{m-1-i}e_i'
=
\sum_{i=0}^{m-1}\sum_{l=1}^{m-i}  (qz_j-1)^ip^{m-1-i}(-1)^{l-1}z_j^{-l}e_{i+l}
\end{aligned}
\]
and
\[
\begin{aligned}
\Delta^T_j
&=
\prod_{k:k\neq j}(qz_jz_k-z_j+p)\\
&=
\sum_{i=0}^{m-1} (qz_j)^i(-z_j+p)^{m-1-i}e_i'
=
\sum_{i=0}^{m-1}\sum_{l=0}^{i}  (qz_j)^i(-z_j+p)^{m-1-i}(-1)^lz_j^le_{i-l}.
\end{aligned}
\]

Further, by exchanging indices and using the formula $1+X+X^2+\dots+X^N=\frac{1-X^{N+1}}{1-X}$, we can significantly simplify the above expressions as follows:
\begin{equation}\label{eq:Delta_Delta_T}
\begin{aligned}
\Delta_j&=
\frac{p^{m}}{(z_j-1)(qz_j-p)}\sum_{l=0}^m e_l\left\{
p^{-l}(qz_j-1)^l
-
(-z_j^{-1})^{l}
\right\},\\
\Delta_j^T&=\frac{q^mz_j^m}{(z_j-1)(qz_j-p)}\sum_{l=0}^me_{m-l} 
\left\{
q^{-l}(pz_j^{-1}-1)^l-(-z_j)^l
\right\}.
\end{aligned}
\end{equation}
Although these expressions take the form of rational functions, they are in fact polynomials in $z_j^{\pm 1}$ with coefficients in $\ZZ[e_1,\dots,e_m,p]$.

\begin{rem}
By using the trivial equations 
\[
\sum_{l=0}^me_l(-z_j^{-1})^l=
\sum_{l=0}^me_{m-l}(-z_j)^l=0,
\]
\eqref{eq:Delta_Delta_T} can be further simplified as
\begin{equation}\label{eq:Delta_simplified}
\begin{aligned}
\Delta_j&=
\frac{p^{m}}{(z_j-1)(qz_j-p)}\sum_{l=0}^m e_l\left(
\frac{qz_j-1}{p}
\right)^l,\\
\Delta_j^T&=\frac{q^mz_j^m}{(z_j-1)(qz_j-p)}\sum_{l=0}^me_{m-l} 
\left(
\frac{pz_j^{-1}-1}{q}
\right)^l.
\end{aligned}
\end{equation}
Here, substituting $e_i = \sum_{1\leq i_1<\dots<i_k\leq m} \prod_{j=1}^k z_{i_j}$ transforms \eqref{eq:Delta_simplified} into a polynomial in $z_1,\dots,z_m$. However, since \eqref{eq:Delta_Delta_T} yields a polynomial without such a substitution, we often prefer using \eqref{eq:Delta_Delta_T} rather than \eqref{eq:Delta_simplified}.
\end{rem}

Let $\varphi(\lambda;p)$ and $\psi(\lambda;p)$ be the polynomials in $\lambda$ defined by
\[
\begin{aligned}
\varphi(\lambda;p)&=
\lambda^n\cdot 
\frac{p^{m}}{(\lambda-1)(q\lambda-p)}\sum_{l=0}^m e_l\left\{
p^{-l}(q\lambda-1)^l
-
(-\lambda^{-1})^{l}
\right\},\\
\psi(\lambda;p)&=
\frac{q^m\lambda^m}{(\lambda-1)(q\lambda-p)}\sum_{l=0}^me_{m-l} 
\left\{
q^{-l}(p\lambda^{-1}-1)^l-(-\lambda)^l
\right\}.
\end{aligned}
\]
Then, the polynomial form of the ASEP-Bethe equations \eqref{eq:ASEP_poly} can be rewritten as 
\begin{equation}\label{eq:relation_varphi_0}
\varphi(z_j;p)+(-1)^mg\psi(z_j;p)=0,\qquad (j=1,\dots,m).
\end{equation}

Let $u_{i}(e_1,\dots,e_m;p)$ and $v_{i}(e_1,\dots,e_m;p)$ be the coefficients of the expansion
\[
\begin{aligned}
\varphi(\lambda;p)&=
(-1)^{m+n-2}
\sum_{i=n-m}^{m+n-2}(-1)^{i}u_{i}(e_1,\dots,e_m;p)\lambda^i,\\
\psi(\lambda;p)&=
(-1)^{m+n-2}
\sum_{i=0}^{2m-2}
(-1)^{i}
v_{i}(e_1,\dots,e_m;p)\lambda^i.
\end{aligned}
\]
Then, the associated projectivized system \eqref{eq:modif_proj} has the explicit form
\begin{equation}\label{eq:vector_form_ASEP}
\begin{bmatrix}
a_0b_1\\
a_0b_2+a_1b_1\\
\vdots\\
a_{m-1}b_{n-1}+a_{m}b_{n-2}\\
a_{m}b_{n-1}\\
\end{bmatrix}
=
b_0
\begin{bmatrix}
\mathcal{U}_{m+n-2}(a;p)\\
\mathcal{U}_{m+n-3}(a;p)\\
\vdots\\
\mathcal{U}_{n-m}(a;p)\\
\mathbf{0}_{n-m}
\end{bmatrix}
+(-1)^mb_0g
\begin{bmatrix}
\mathbf{0}_{n-m}\\
\mathcal{V}_{2m-2}(a;p)\\
\mathcal{V}_{2m-3}(a;p)\\
\vdots\\
\mathcal{V}_{0}(a;p)
\end{bmatrix}.
\end{equation}
Here, $
([a_0:\dots :a_m],[b_0:\dots:b_{n-1}],p)\in \PP^m\times \PP^{n-1}\times \CC$ is the standard coordinates and
\[
\mathcal{U}_k(a;p)=a_0\cdot u_k(a_1/a_0,\dots,a_m/a_0;p),\quad
\mathcal{V}_k(a;p)=a_0\cdot v_k(a_1/a_0,\dots,a_m/a_0;p)
\]
are the homogenization of $u_k$ and $v_k$, respectively.
(Note that the degree of $u_k$ and $v_k$ are $1$.)

Let $t$ be a formal variable.
By introducing the generating functions
\[
\begin{aligned}
a(t)&=a_0-a_1t+a_2t^2-\cdots+(-1)^ma_mt^m,\\
b(t)&=b_1-b_2t+b_3t^2-\cdots+(-1)^{n-2}b_{n-1}t^{n-2},
\end{aligned}
\]
the system \eqref{eq:vector_form_ASEP} can be rewritten as
\begin{equation}\label{eq:system_a_and_b_m}
\begin{aligned}
a(t)b(t)&=t^{m+n-2}b_0\left(\varphi(t^{-1};p)+(-1)^mg\psi(t^{-1};p)\right)\\
&=
\frac{b_0}{(1-t)(q-pt)}
\Bigg[
p^{m}t^m\sum_{l=0}^m a_l\left\{
p^{-l}(qt^{-1}-1)^l
-
(-t)^{l}
\right\}\\
&\hspace{12em}+(-1)^mgq^mt^{n}
\sum_{l=0}^ma_{m-l} 
\left\{
q^{-l}(pt-1)^l-(-t^{-1})^l
\right\}
\Bigg].
\end{aligned}
\end{equation}
By multiplying $(1-t)(q-pt)$ and comparing the coefficients of $t^0$ and $t^{m+n}$ on the both sides, we obtain
\begin{equation}\label{eq:a_0_b_1}
a_0b_1=q^{m-1}a_mb_0\qquad\text{ and }\qquad
a_mb_{n-1}=(-1)^np^{m-1}ga_0b_0,
\end{equation}
respectively.

If a pair $(a(t),b(t))$ satisfies \eqref{eq:system_a_and_b_m} with a specified complex number $p$, we often use the notations
\[
a(t;p)=a(t),\qquad b(t;p)=b(t)
\]
to specify the dependence on $p$.

\begin{rem}
If $a(t;p)$ is factorized as
\[
a(t;p)=c(1-y_1t)(1-y_2t)\cdots (1-y_mt),\qquad (c\neq 0)
\]
for some functions $y_i=y_i(p)$ in $p$, the right hand side of \eqref{eq:system_a_and_b_m} is simplified as
\begin{equation}\label{eq:simplified}
\begin{aligned}
&\frac{b_0c}{(1-t)(q-pt)}\Bigg[
\prod_{i=1}^m(pt+y_iq-y_it)
+gt^n
\prod_{i=1}^m(-pt+1-y_iq)\\
&\hspace{15em}-(p^mt^m+gq^mt^{n-m})\prod_{i=1}^m(1-y_it)
\Bigg].
\end{aligned}
\end{equation}
This expression is useful for checking whether a given function is a solution to \eqref{eq:system_a_and_b_m}.
\end{rem}

\subsection{The scheme $S$ of quasi-simple solutions}
\label{sec:def_of_S}

Fix a non-zero complex number $g\neq 0$.
From the $k$-th component of \eqref{eq:vector_form_ASEP}, we define the polynomial $f_k(a,b;p)\in \CC[a_i,b_i,p]$ by
\begin{equation}\label{eq:def_of_f_k}
f_k(a,b;p):=\left(\sum_{\substack{i+j=k,\\ i\geq 0,j\geq 1}}a_ib_j\right) -b_0
\left(\mathcal{U}_{m+n-1-k}(a;p)
+
(-1)^m
g\mathcal{V}_{m+n-1-k}(a;p)
\right),
\end{equation}
which is bi-homogeneous in $a_i$ and $b_i$.

Let 
\[
(a,b,p)=([a_0:\dots:a_m],[b_0:\dots:b_{n-1}],p)\qquad \in \PP^m\times \PP^{n-1}\times \CC
\]
be the standard coordinate of $\PP^m\times \PP^{n-1}\times \CC$.
Let $D_k\subset \PP^m\times \PP^{n-1}\times \CC$ denote the Cartier divisor defined by $f_k$.
Then, the scheme theoretic intersection 
\[
S:=D_1\cap D_2\cap \dots \cap D_{m+n-1}
\]
can be seen as a parameter space of all quasi-simple solutions to \eqref{eq:ASEP_poly}.
Since $\PP^m\times \PP^{n-1}\times \CC$ is $(m+n)$-dimensional, the expected dimension of $S$ is $1=(m+n)-(m+n-1)$.  

Hereafter, we study the geometric properties of $S$.
For this, we focus on intersections of $S$ and hyperplanes, and their neighborhoods.
For brevity, we will use the shorthand notations
\[
\begin{aligned}
\{a_i=c\}&=\{(a,b,p)\in  \PP^m\times \PP^{n-1}\times\CC: a_i=c \},\\
\{b_i=c\}&=\{(a,b,p)\in  \PP^m\times \PP^{n-1}\times\CC: b_i=c \},\\
\{p=c\}&=\{(a,b,p)\in  \PP^m\times \PP^{n-1}\times\CC: p=c \}
\end{aligned}
\]
for $c\in \CC$.

\begin{lemma}\label{lemma:b_0_is_not_0}
We have $S\cap \{b_0=0\}=\emptyset$.
\end{lemma}
\begin{proof}
If there were a point 
\[
([a_0:\dots:a_m],[0:b_1:\dots:b_{n-1}],p)\qquad \in S\cap \{b_0=0\},
\] 
it would follow from \eqref{eq:system_a_and_b_m} that $a(t;p)b(t;p)=0$.
However, this is impossible because neither $a_0=\dots=a_n=0$ nor $b_0=\dots=b_{n-1}=0$ can hold.
\end{proof}

By Lemma \ref{lemma:b_0_is_not_0}, we can assume
\[
b_0=1
\]
whenever we deal with $S$.

\subsection{Irreducible components of $S$}\label{sec:Irr_of_S}

The polynomial $\Delta(z_1,\dots,z_m)$ (Eq.~\eqref{eq:def_of_Delta}) defines the locus on which the ASEP-Bethe equations \eqref{eq:ASEP_Bethe} are undefined.
Since $\Delta(z_1,\dots,z_m)$ is a symmetric polynomial in $z_1,\dots,z_m$, there exists a polynomial $E(e_1,\dots,e_n;p)\in \ZZ[e_1,\dots,e_n,p]$ satisfying
\begin{equation}\label{eq:def_of_E}
E(e_1,\dots,e_n;p)=\Delta(z_1,\dots,z_n)=\prod_{(j,k):j\neq k}(qz_jz_k-z_k+p).
\end{equation}
Let 
\[
\mathcal{E}(a_0,\dots,a_m;p):=a_0^{\deg(E)}\cdot E(a_1/a_0,\dots,a_m/a_0;p)
\]
be the homogenization of $E$.

Although finding the explicit form of $\mathcal{E}$ is not straightforward, it is not needed in the following discussion. 
Everything we need to know about the set $\{\mathcal{E}=0\}$ is summarized in the following proposition.   

\begin{prop}\label{prop:E_set}
The algebraic set $S\cap \{\mathcal{E}=0\}$ is contained in
\begin{align*}
&\{a_0=0\}\cup \{a_m=0\}\cup \{a_0-a_1+\dots+(-1)^ma_m=0\}\\
&\qquad \cup 
\{a_0p^m-a_1p^{m-1}q+\dots+(-1)^ma_mq^m=0\}\cup \bigcup_{l=2}^m\{(p/q)^l=1\}. 
\end{align*}
\end{prop}
\begin{proof}
Assume a point $([1:a_1:\dots:a_m],[1:b_1:\dots:b_{n-1}],p)$ is contained in 
\[
S\cap \{\mathcal{E}=0\}\cap \{a_0\neq 0\}\cap \{a_m\neq 0\}.
\]
Let $z_1,\dots,z_m$ be the roots of $z^m-a_1z^{m-1}+\dots+(-1)^ma_m$.
Then, we have $\Delta(z_1,\dots,z_n)=0$ and $z_1z_2\dots z_m\neq 0$ by assumption.
By \eqref{eq:def_of_E}, there exists at least one pair $(k,j)$ satisfying $F(z_k,z_j)=qz_kz_j-z_j+p=0$.
Moreover, from \eqref{eq:ASEP_poly}, there must exist an index $i\neq k$ such that $F(z_i,z_k)=0$.
Repeating this procedure, we can find a sequence $j_1,j_2,\dots,j_l$ such that
\begin{equation}\label{eq:F_is_0}
F(z_{i_1},z_{i_2})=F(z_{i_2},z_{i_3})=\dots=F(z_{i_l},z_{i_1})=0.
\end{equation}
If $p=0$ or $p=1$, we have $z_{i_1}=z_{i_2}=\dots=1$ from \eqref{eq:F_is_0}, which implies $a_0-a_1+\dots+(-1)^ma_m=0$.

Let $p\notin \{0,1\}$. 
We can assume $i_k=k$ without loss of generality.
For the linear fractional transformation $\varphi$ defined by $\varphi(t)=\frac{p}{1-qt}$, it follows from \eqref{eq:F_is_0} that 
\[
z_{2}=\varphi(z_{1}),\quad
z_{3}=\varphi(z_{2}),\quad
\dots,\quad
z_{1}=\varphi(z_{l}).
\]
Then $z_1$ is a fixed point of $\varphi^l=\overbrace{\varphi\circ \dots \circ \varphi}^{l}$.

If the matrix 
$\begin{bmatrix}
0&p\\-q&1
\end{bmatrix}^l$ is not a scalar multiple of $\mathrm{Id}=\begin{bmatrix}
1&0\\0&1
\end{bmatrix}$, $z_1$ is also a fixed point of $\varphi$; hence, we have $z_1\in \{1,p/q\}$, which implies \[
(a_0-a_1+\dots+(-1)^ma_m)(a_0p^m-a_1p^{m-1}q+\dots+(-1)^ma_mq^m)=0.
\]
If $\begin{bmatrix}
0&p\\-q&1
\end{bmatrix}^l$ is a scalar multiple of $\mathrm{Id}$, then we have $(p/q)^l=1$.
\end{proof}

In general, the scheme $S$ consists of finite number of irreducible components.
Since $S$ is defined by $(m+n-1)$ effective Cartier divisors, each component of $S$ is of dimension $\geq 1$.
By \eqref{eq:e_m_is_const}, each irreducible component of $S$ must be contained in at least one of the three hypersurfaces:
\[
\{a_m^n-g^ma_0^n=0\},\quad \{a_0=0\},\quad \{\mathcal{E}=0\}.
\]

\begin{defi}
Let $Z$ be an irreducible component of $S$.
\begin{itemize}
\item $Z$ is called an infinite component if $Z$ is contained in $\{a_0=0\}$.
\item $Z$ is called an indefinite component if $Z$ is contained in $\{\mathcal{E}=0\}$.
\item $Z$ is called a vertical component if $Z$ is contained in $\{p=c\}$ for some $c\in \CC$.
\item $Z$ is called an essential component if $Z$ is not an infinite, indefinite, nor vertical component.
All essential components are contained in $\{a_m^n-g^ma_0^n=0\}$.
\end{itemize}
\end{defi}
As we will see in \S \ref{sec:special_fiber}, the four types of irreducible components can be completely classified by examining how they intersect with the hyperplane $\{p=0\}$.
Note that, for an irreducible component $Z\subset S$ that is not vertical, it follows from \eqref{eq:a_0_b_1} that
\begin{equation}\label{eq:a_0_equiv_a_m}
Z\subset \{a_0=0\}\iff Z\subset \{a_m=0\}.
\end{equation}

The following proposition deals with the existence of the steady state solution $z_1=z_2=\dots=z_m=1$.
As mentioned in \S \ref{sec:ASEP_projec}, the steady state solution is the sole exception and requires careful handling.

\begin{prop}\label{prop:steady}
If $g=1$, \eqref{eq:system_a_and_b_m} admits the special solution
\begin{equation}\label{eq:steady_state}
a(t;p)=a_0(1-t)^m,\qquad 
b(t;p)=
\frac{q^m-p^mt^m}{q-pt}
\frac{1-t^{n-m}}{1-t}
,\qquad (a_0\neq 0).
\end{equation}
(Note that $b(t;p)$ is in fact a polynomial in $t$ of degree $n-2$.)
Conversely, if \eqref{eq:system_a_and_b_m} admits a solution $(a(t;p),b(t;p))$ satisfying $a(t;p)=a_0(1-t)^m$ for all $p\in \CC$, then we have $g=1$.
\end{prop}
\begin{proof}
We can apply the normalization $a_0=b_0=1$ here.
Substituting $a(t;p)=(1-t)^m$ and $y_1=\dots=y_m=1$ to \eqref{eq:system_a_and_b_m} and \eqref{eq:simplified}, we have 
\[
(1-t)^mb(t;p)=\frac{(q^m-p^mt^m)(1-gt^{n-m})(1-t)^m}{(1-t)(q-pt)}.
\]
Hence, we have
\[
b(t;p)=\frac{q^m-p^mt^m}{q-pt}\cdot \frac{1-gt^{n-m}}{1-t}.
\]
If $g=1$, $b(t;p)$ is a polynomial; hence, this is a solution.
Conversely, if $b(t;p)$ is a polynomial for arbitrary $p$, then $g$ must satisfy $g=1$.
\end{proof}

The function $a(t;p)=a_0(1-t)^m$ corresponds to the Bethe solution $z_1=z_2=\dots=z_m=1$.
Let $Z_{\mathrm{steady}}\subset \PP^m\times \PP^{n-1}\times \CC$ be the subvariety defined by
\[
Z_{\mathrm{steady}}:=\{([a_0:\dots:a_m],[b_0,\dots,b_{n-1}],p):
\text{The generating functions $a(t)$ and $b(t)$ satisfy \eqref{eq:steady_state}.}
\}
\]
We call $Z_{\mathrm{steady}}$ the \textit{component of steady state}.

\begin{rem}
Regardless of the value of $g$, the polynomial form of the ASEP-Bethe equations \eqref{eq:ASEP_poly} admits the solution $z_1=\dots=z_m=1$.
However, Proposition \ref{prop:steady} claims that if the system admits a quasi-simple solution $z_1=\dots=z_m=1$ for all $p\in \CC$, then $g$ must satisfy $g=1$.
This fact precisely reflects the requirement of physics that the quasi-periodic ASEP admits the steady state only if it is a properly periodic ($g=1$).
\end{rem}

\subsection{Special fiber $S|_0$}\label{sec:special_fiber}

Let 
\begin{equation}\label{eq:projection_def}
\pi:\PP^{m}\times \PP^{n-1}\times \CC\to \CC
\end{equation}
be the projection to the third component.
For a closed subscheme $A$ of $\PP^{m}\times \PP^{n-1}\times \CC$ and $p\in \CC$, we write 
\[
A|_p:=A\cap \pi^{-1}(p).
\]

We consider the special fiber $S|_0$, which corresponds to the TASEP case \eqref{eq:B_P_TASEP}.
To understand the global structure of $S$, it is necessary to analyze the local behavior in a neighborhood of $p = 0$, in addition to a pointwise analysis at $p = 0$. 

\subsubsection{Analysis at $p=0$}

Substituting $p=0$, $q=1$ and $b_0=1$ to 
\eqref{eq:system_a_and_b_m} yields
\[
\begin{aligned}
a(t;0)b(t;0)
&
=\frac{1}{1-t}
\left[a_m(1-t)^{m}
+(-1)^m
gt^{n}
\sum_{l=0}^ma_{m-l} 
\left\{
(-1)^l-(-t^{-1})^l
\right\}
\right]\\
&=
\frac{1}{1-t}
\left[a_m(1-t)^{m}
+
g\{
a(1;0)t^n-a(t;0)t^{n-m}
\}
\right].
\end{aligned}
\]
By putting 
\begin{equation}\label{eq:def_of_h}
h(t;0):=b(t;0)(1-t)+gt^{n-m}, 
\end{equation}
we obtain the simple equation 
\begin{equation}\label{eq:a_and_h}
a(t;0)h(t;0)=a_m(1-t)^m+g\cdot a(1;0)t^n.
\end{equation}

The special fiber $S|_0$ consists of all points whose standard coordinates
\begin{equation}\label{eq:coordinate}
([a_0:\dots:a_m],[1:b_1:\dots:b_{n-1}],0)
\end{equation}
satisfies \eqref{eq:a_and_h}.
\begin{prop}\label{prop:bound_points}
The homogeneous coordinates \eqref{eq:coordinate} satisfies exactly one of the following (i) and (ii):
\begin{itemize}
\item[(i)] $a_m\neq 0$ and $b_{n-1}=b_{n-2}=\dots=b_{n-m+1}=0$,
\item[(ii)] for some $1\leq k\leq m-1$,
\begin{equation}\label{eq:homo_corr}
a_i=
\begin{cases}
1 & (i=m-k),\\
0 & (i\neq m-k),
\end{cases}\qquad
b_i=
\begin{cases}
(-1)^{i}g & (n-m+1\leq i\leq n-m+k),\\
0 & (\text{otherwise}).
\end{cases}
\end{equation}
\end{itemize}
\end{prop}
\begin{proof}
Let $h(t;0)=h_0-h_1t+\dots+(-1)^{n-1}h_{n-1}t^{n-1}$.
Then, by comparing the degrees of the polynomials on both sides of \eqref{eq:a_and_h}, we find that there exists $0\leq k\leq m-1$ such that
\[
a_m=a_{m-1}=\dots=a_{m-k+1}=0,\quad a_{m-k}\neq 0,\quad \text{and}\quad  h_{n-1}=\dots=h_{n-m+k+1}=0.
\]
If $k=0$, (i) follows from \eqref{eq:def_of_h}.
If $k>0$, we obtain $a_m=0$ and 
\[
a(t;0)h(t;0)=g\cdot a(1;0)t^n
\]
from \eqref{eq:a_and_h}.
As $\deg(a(t;0))=m-k$, we have
$a(t;0)=(-1)^{m-k}a_{m-k}t^{m-k}$ and $h(t;0)=gt^{n-m+k}$,
which implies
\[
b(t;0)
=
g\frac{t^{n-m+k}-t^{n-m}}{1-t}=-gt^{n-m}(1+t+\dots+t^{k-1}).
\]
The equations \eqref{eq:homo_corr} are obtained by taking the coefficients of the polynomials $a(t;0)$ and $b(t;0)$.
\end{proof}

There are some distinguished points contained in $S|_0$.
By Proposition \ref{prop:bound_points} (ii) and \eqref{eq:homo_corr}, $S|_0$ contains $(m-1)$ points $P_1,P_2,\dots,P_{m-1}$, where the homogeneous coordinates of $P_k$ are given by
\begin{equation}\label{eq:def_of_P_k}
\begin{aligned}
&([a_0:\dots:a_m],[b_0:b_1:\dots:b_{n-1}],p)\\
&=
\left([0:\dots:\overbrace{1}^{a_{m-k}}:\dots:0],
[1:\overbrace{0:\dots:0}^{n-m}:\overbrace{(-1)^{n-m+1}g:\dots:(-1)^{n-m+k}g}^k:0:\dots:0],0
\right).
\end{aligned}
\end{equation}
In addition, if $g=1$, $S|_0$ contains the point $P_{\mathrm{steady}}$ whose homogeneous coordinates are
\begin{equation}\label{eq:def_of_P_steady}
\begin{aligned}
&([a_0:\dots:a_m],[b_0:b_1:\dots:b_{n-1}],p)\\
&=
\left(\left[1:\binom{m}{1}:\binom{m}{2}:\dots:1\right],
[1:1:-1:\dots:(-1)^{n-m-1}:\overbrace{0:\dots:0}^{m-1}],0
\right),
\end{aligned}
\end{equation}
which is derived from \eqref{eq:steady_state} by putting $p=0$.

\begin{lemma}\label{lemma:a_1_0}
Assume that a pair $(a(t;0),b(t;0))$ satisfies \eqref{eq:def_of_h} and \eqref{eq:a_and_h}.
Then, we have the following (i) and (ii).
\begin{enumerate}
\item[(i)] If $g\neq 1$, then $a(1;0)\neq 0$.
\item[(ii)] If $g=1$, then $a(1;0)=0$ implies 
\[
\left(a(t;0),b(t;0)\right)=\left(
a_m(1-t)^m,1+t+\dots+t^{n-m-1}
\right),
\]
which corresponds to the point $P_{\mathrm{steady}}$ (Eq.\eqref{eq:def_of_P_steady}).
\end{enumerate}
\end{lemma}
\begin{proof}
Suppose $a(1;0)=0$.
Then, it follows from \eqref{eq:a_and_h} that
\begin{equation}\label{eq:sub_eq}
a(t;0)h(t;0)=a_m(1-t)^m.
\end{equation}
If $a_m=0$, \eqref{eq:sub_eq} leads $h(t;0)=0$ because $a_0,a_1,\dots,a_m$ cannot vanish simultaneously.
Then, we have $b(t;0)(1-t)=-gt^{n-m}$, which contradicts $g\neq 0$.
If $a_m\neq 0$, it follows from \eqref{eq:sub_eq} that $a(t;0)=a_m(1-t)^m$ and $h(t;0)=1$.
From \eqref{eq:def_of_h}, we have $g=1$ and $b(t;0)=1+t+\dots+t^{n-m-1}$.
\end{proof}

While it is easy to expect that the dimension of $S|_0$ is $0$, proving that this actually holds is non-trivial.
From Proposition \ref{prop:bound_points}, we have
\begin{equation}\label{eq:intersection_S_0_inft}
S|_0\cap \{a_0=0\}=\{P_1,\dots,P_{m-1}\}.
\end{equation}

\begin{theorem}\label{thm:0dim}
$S|_0$ is $0$-dimensional.
\end{theorem}
\begin{proof}
By Lemma \ref{lemma:a_1_0} and \eqref{eq:intersection_S_0_inft}, it suffices to show that 
\[
\mathcal{Q}:=S|_0\cap\{a_m\neq 0\}\cap \{a_0-a_1+\dots+(-1)^ma_m\neq 0\}
\]
is a finite set.
We can assume $a_m=(-1)^m$ by normalization without loss of generality.
Let $\mathcal{X}(g)$ and $\mathcal{Y}(g)$ be the two sets 
\[
\begin{gathered}
\mathcal{X}(g)=\left\{(a(t),h(t)):
\begin{aligned}
&a(t)h(t)=(t-1)^m+g\cdot a(1)t^n,\\ &\deg(a(t))=m,\quad \text{$a(t)$ is monic},\quad a(1)\neq 0
\end{aligned}
\right\},\\
\mathcal{Y}(g)=
\left\{(a(t),h(t),C):
\begin{aligned}
&a(t)h(t)=(t-1)^m+Ct^n,\\
&\deg(a(t))=m,\quad \text{$a(t)$ is monic},\quad h(1)=g,\ C\neq 0 
\end{aligned}
\right\}.
\end{gathered}
\]
By \eqref{eq:def_of_h} and \eqref{eq:a_and_h}, there exists a one-to-one correspondence between $\mathcal{X}(g)$ and $\mathcal{Q}$.
Since the map $\mathcal{Y}(g)\to \mathcal{X}(g)$ sending $(a(t),h(t),C)$ to $(a(t),h(t))$ is bijective, it suffices to prove $\mathcal{Y}(g)$ to be finite.

Let 
\begin{equation}\label{eq:factor_on}
(t-1)^m+Ct^{n}=C\prod_{i=1}^n(t-\alpha_i)
\end{equation}
be the factorization, where $\alpha_1=\alpha_1(C),\dots,\alpha_n=\alpha_n(C)$ are the roots of $(t-1)^m+Ct^{n}$.
Consider the projection $j:\mathcal{Y}(g)\to \CC$ that sends $(a(t),h(t),C)$ to $C$.
For any element of the fiber $j^{-1}(C)$, there exists at least one subset $J\subset \{1,\dots,n\}$ with $|J|=n-m$ satisfying
\begin{equation}\label{eq:valuation}
h(1)=C\prod_{j\in J}(1-\alpha_j)=g.
\end{equation}
Let $\mathcal{C}(g)\subset \CC\setminus \{0\}$ be the subset of complex numbers $C$ such that $j^{-1}(C)$ is nonempty.
Since each branch of the multi-valued map
\begin{equation}\label{eq:decomposition_C}
C\mapsto \left\{C\prod_{j\in J}(1-\alpha_j(C)):|J|=n-m\right\}
\end{equation}
is an algebraic function\footnote{The image of the multi-valued map \eqref{eq:decomposition_C} consists of the root of 
\begin{equation*}
P(T;C):=\prod_{|J|=n-m}\left(T-C\prod_{j\in J}(1-\alpha_j)\right),
\end{equation*}
which is a polynomial in $T$ of degree $\binom{n}{n-m}$.
Since all coefficients of $P(T;C)$ is symmetric in $\alpha_1,\dots,\alpha_n$, $P(T;C)$ is in fact an element of $\ZZ[C^{\pm 1},T]$.
The structure of the Riemann surface (= algebraic curve) defined by $P(T;C)=0$ is deeply studied in the previous works by Prolhac~\cite{prolhac2020riemann,prolhac2021riemann} and Iwao-Motegi~\cite{iwao2025bethe}.
} 
in $C$, the set $\mathcal{C}(g)$ must satisfy (i) $|\mathcal{C}(g)|<+\infty$ or (ii) $\mathcal{C}(g)=\CC\setminus \{0\}$.
However, Case (ii) cannot occur because, if we let $C \to +\infty$ continuously, the root $\alpha_i$ would converge to $0$ as $\alpha_i = O(C^{-1/n})$, which contradicts $C\prod_{j\in J}(1-\alpha_j)=g$.
Therefore, the theorem follows from the inequality $|\mathcal{Y}(g)| \leq |\mathcal{C}(g)|\times \binom{n}{n-m}<+\infty$, which holds by \eqref{eq:valuation}.
\end{proof}

From Theorem \ref{thm:0dim}, we obtain the following important fact:

\begin{cor}
An irreducible component $Z\subset S$ that intersects the hyperplane $\{p=0\}$ is $1$-dimensional.
\end{cor}
\begin{proof}
As mentioned above, we have $\dim Z\geq 1$.
On the other hand, from Theorem \ref{thm:0dim}, it follows that $\dim Z\leq \mathrm{codim}(\{p=0\})+\dim (Z|_0)=1+0=1$, which concludes $\dim Z=1$.
\end{proof}

\subsubsection{Analysis in a neighborhood of $p=0$}\label{sec:aroung_p_0}

It is straightforward to check that $P_k$ is contained in $\{a_m^n-g^ma_0^n=0\}$.
However, this does not mean that $P_k$ is a point in an essential component.
In fact, if $P_k$ were contained in an essential component $Z$, it would follow from \eqref{eq:a_0_b_1} that
\[
0=a_0b_1-q^{m-1}a_m=a_0(b_1-cq^{m-1}),\qquad \text{ on $Z$}
\]
for some $c\neq 0$ satisfying $c^n=g^m$.
Since $Z\not\subset\{a_0=0\}$, $a_0$ is not identically $0$ on $Z$; hence, $b_1$ must satisfy $b_1=cq^{m-1}$. 
However, this contradicts $\lim\limits_{p\to 0} b_1=0$ (see \eqref{eq:def_of_P_k}).
Therefore, $P_k$ is indeed not a point in an essential component.


As this example shows, to determine which points belong to essential components, we need to examine the behavior of $S$ in a neighborhood of $p = 0$ as well.
In particular, we focus on the intersections of indefinite components and $\{p=0\}$.
The following Proposition \ref{lemma:S_ex_0} indicates that the possible cases of such intersections are limited.





\begin{prop}\label{lemma:S_ex_0}
For an indefinite component $Z\subset S$, we have
\[
\begin{cases}
Z|_0\subset \{P_1,\dots,P_{m-1}\} & (g\neq 1),\\
Z|_0\subset \{P_1,\dots,P_{m-1}\}\cup \{P_{\mathrm{steady}}\} & (g=1).
\end{cases}
\]
\end{prop}

\begin{proof}
Let $([a_0:\dots:a_m],[1:b_1:\dots:b_{n-1}])\in Z|_0$.
As $Z$ is indefinite, it follows from Proposition \ref{prop:E_set} and \eqref{eq:a_0_equiv_a_m} that $a_m=0$ or $a_0-a_1+\dots+(-1)^ma_m=0$.
Therefore, it holds that
\[
a_m\cdot a(1;0)=
a_m(a_0-a_1+\dots+(-1)^ma_m)=0.
\]
From Lemma \ref{lemma:a_1_0}, we have 
\begin{itemize}
\item [(i)] $g\neq 1\Rightarrow a_m=0$, and
\item [(ii)] $g=1\Rightarrow \begin{cases}a_m=0, or \\(a(t;0),b(t;0))\text{ corresponds to $P_{\mathrm{steady}}$}.\end{cases}$
\end{itemize}
Hence, we have
\[
Z|_0\subset
\begin{cases}
\{a_m=0\} & (g\neq 1),\\
\{a_m=0\}\cup \{P_{\mathrm{steady}}\} & (g= 1).
\end{cases}
\]
From $Z|_0\cap \{a_m=0\}\subset \{P_1,\dots,P_{m-1}\}$, which follows from Proposition \ref{prop:bound_points}, we conclude the proposition.
\end{proof}

Among indefinite components, $Z_{\mathrm{steady}}$ is solely special.
Set $g=1$.
Let $Z'$ be an arbitrary irreducible component containing $P_{\mathrm{steady}}$.
Since $Z'$ is $1$-dimensional and does not contained in $\{p=0\}$, $Z'$ is locally expressed as a graph of local coordinate functions
\[
A_k(p):=\frac{a_k}{a_0}(p),\qquad B_k(p):=b_k(p),
\] 
where $A^{(l)}_k(p)$ and $B_k(p)$ are Puiseux series
\[
A_k(p)=\binom{m}{k}+\alpha'_{k}p^\nu+\alpha''_{k}p^{2\nu}+\dots,\quad
B_k(p)=\beta_{k}+\beta'_{k}p^\nu+\beta''_{k}p^{2\nu}+\dots
\]
for some $\nu\in \mathbb{Q}_{>0}$.

Define the generating functions
\[
A(t;p)=1+\sum_{k=1}^m(-1)^kA_k(p)t^k,\quad B(t;p)= \sum_{k=0}^{n-2}(-1)^kB_{k+1}(p)t^k. 
\]
Then, the pair $(A(t;p),B(t;p))$ satisfies \eqref{eq:system_a_and_b_m} for sufficiently small $|p|$.
Let 
\[
\delta A(t;p):=A(t;p)-(1-t)^m,\qquad
\delta B(t;p):=B(t;p)-\frac{q^m-p^mt^m}{q-pt}
\frac{1-t^{n-m}}{1-t}
\]
be the difference of $(A(t;p),B(t;p))$ and the steady state solution.
Expand these functions as
\begin{equation}\label{eq:expansion_delta}
\begin{aligned}
\delta A(t;p)
&=A^{(1)}(t)p^{\nu}+A^{(2)}(t)p^{2\nu}+\cdots,\\ 
\delta B(t;p)
&=B^{(1)}(t)p^{\nu}+B^{(2)}(t)p^{2\nu}+\cdots.
\end{aligned}
\end{equation}
By definition of $\delta A$, we have $\delta A(0;p)=0$.

\begin{prop}\label{prop:Z_steady_unique}
$Z_{\mathrm{steady}}$ is the only irreducible component of $S$ that intersects $\{p=0\}$ at $P_{\mathrm{steady}}$.
\end{prop}
\begin{proof}
We prove the proposition by contradiction.
Assume $P_{\mathrm{steady}}$ were contained in some irreducible component $Z'\neq Z_{\mathrm{steady}}$.
Then, there exists some $l\in \ZZ_{>0}$ such that at least one of the polynomials $A^{(l)}(t)$ and $B^{(l)}(t)$ in \eqref{eq:expansion_delta} is not $0$.
Let $l$ be the least number of such integers.
Let
\[
A^{(l)}(t):=a_0'-a_1't+a_2't-\dots+a'_mt^m.
\]
Since $\delta A(0;p)=0$, we have $A^{(l)}(0)=0$.

By substituting $A(t;p)=(1-t)^m+\delta A(t;p)$ and $B(t;p)=\frac{q^m-p^mt^m}{q-pt}
\frac{1-t^{n-m}}{1-t}+\delta B(t;p)$ to the left hand side of \eqref{eq:system_a_and_b_m}, we obtain
\[
(1-t)^m\frac{q^m-p^mt^m}{q-pt}\frac{1-t^{n-m}}{1-t}+
\left(
A^{(l)}(t)\frac{q^m-p^mt^m}{q-pt}\frac{1-t^{n-m}}{1-t}+(1-t)^mB^{(l)}(t)
\right)p^{l\nu}+o(p^{l\nu}).
\]
On the other hand, since the right hand side of \eqref{eq:system_a_and_b_m} is a linear function in $a_0,a_1,\dots,a_m$, we obtain
\[
(\text{RHS of \eqref{eq:system_a_and_b_m}})|_{a_k\mapsto \binom{m}{k}}
+
(\text{RHS of \eqref{eq:system_a_and_b_m}})|_{\substack{
a_0\mapsto 0\\
a_k\mapsto a'_{k}}}p^{l\nu}
+
o(p^{l\nu}).
\]
Taking $p\to 0$ and comparing their coefficients of $p^{l\nu}$, we have
\[
A^{(l)}(t)\frac{1-t^{n-m}}{1-t}+(1-t)^mB^{(l)}(t)=
\frac{1}{1-t}
\left[a'_{m}(1-t)^{m}
+
A^{(l)}(1)t^n-A^{(l)}(t)t^{n-m}
\right],
\]
which yields
\begin{equation}\label{eq:alpha_and_beta}
A^{(l)}(t)+(1-t)^{m+1}B^{(l)}(t)
=a'_{m}(1-t)^m+A^{(l)}(1)t^n.
\end{equation}

Differentiating the both sides of \eqref{eq:alpha_and_beta} $m$ times and substituting $t=1$, we obtain $a'_m=a'_m+\binom{n}{m}A^{(l)}(1)$, which implies $A^{(l)}(1)=0$.
Substituting $A^{(l)}(1)=0$ to \eqref{eq:alpha_and_beta}, we obtain
\[
A^{(l)}(t)\equiv a'_m(1-t)^m\mod{(1-t)^{m+1}}.
\]
However, since $A^{(l)}(t)$ is a polynomial in $t$ of degree $\leq m$, we have $A^{(l)}(t)= a'_m(1-t)^m$ and $B^{(l)}(t)=0$.
However, as $A^{(l)}(0)=0$, we must conclude $A^{(l)}(t)=0$, which contradict the definition of $l$.
\end{proof}

By virtue to Propositions \ref{prop:E_set}, \ref{lemma:S_ex_0}, and \ref{prop:Z_steady_unique}, we can completely exclude the effect of indefinite components, except for the steady state component $Z_{\mathrm{steady}}$, from our counting arguments.
Namely, we have the following proposition:
\begin{prop}\label{prop:no_indefinite}
$Z_{\mathrm{steady}}$ is the only indefinite component $Z$ that satisfies $Z\cap \{p=0\}\neq \emptyset$ and $Z\not\subset \{a_0=0\}$.
\end{prop}
\begin{proof}
We prove the claim by contradiction.
Let $Z\subset S$ be an indefinite component with $Z\cap \{p=0\}\neq \emptyset$ and $Z\not\subset \{a_0=0\}$.
Assume $Z\neq Z_{\mathrm{steady}}$.
By Propositions \ref{lemma:S_ex_0} and \ref{prop:Z_steady_unique}, the intersection $Z\cap \{p=0\}$ is a nonempty set contained in $\{P_1,\dots,P_{m-1}\}$.
Therefore, by Proposition \ref{prop:E_set} and \eqref{eq:a_0_equiv_a_m}, $Z$ must be contained in the hypersurface $\{a_0p^m-a_1p^{m-1}q+\dots+(-1)^ma_mq^m=0\}$.
This implies that the equation \eqref{eq:system_a_and_b_m} admits a solution $(a(t;p),b(t;p))$ satisfying $(q-pt)\,|\,a(t;p)$ for all $p$.

Set $\theta:=p/q$.
Let
\[
\begin{aligned}
a(t;p)&=
c(1-\theta t)^k(1-s_1t)(1-s_2t)\cdots(1-s_{m-k}t)\\
&=cq^{-k}(q- pt)^k(1-s_1t)(1-s_2t)\cdots(1-s_{m-k}t)
,
\end{aligned}
\qquad s_i\neq \theta
\]
be the factorization.
Substituting $y_i=s_i$ $(1\leq i\leq m-k)$ and $y_{m-k+1}=\dots=y_m=\theta$ to \eqref{eq:system_a_and_b_m} and \eqref{eq:simplified}, we have 
\[
\begin{aligned}
&(1-\theta t)^k
\left(\prod_{i=1}^{m-k}(1-s_it)\right)\cdot b(t;p)\\
&=\frac{(q-pt)^k}{(1-t)(q-pt)}
\Bigg[
\theta^k
\prod_{i=1}^{m-k}(pt+s_iq-s_it)+gt^n\prod_{i=1}^{m-k}(-pt+1-s_iq)\\
&\hspace{18em}-(p^mt^m+gq^mt^{n-m})q^{-k}\prod_{i=1}^{m-k}(1-s_it)
\Bigg].
\end{aligned}
\]
Multiplying $(1-t)(q-pt)^{1-k}$ on both sides and substituting $t=\theta^{-1}=q/p$, we obtain
\[
\begin{aligned}
0
&=
\theta^k
\prod_{i=1}^{m-k}(q+s_iq-s_i\theta^{-1})+
g\theta^{-n}\prod_{i=1}^{m-k}(-q+1-s_iq)
-(q^m+gq^np^{-n+m})q^{-k}\prod_{i=1}^{m-k}(1-s_i\theta^{-1})\\
&=\prod_{i=1}^{m-k}(p-s_iq)\times
\left[
\theta^{-m+2k}+g
\theta^{-n}
-\theta^{-m+k}-g\theta^{-n+k}
\right]\\
&=\prod_{i=1}^{m-k}(p-s_iq)\times\theta^{-n}(\theta^k-1)(\theta^{n-m+k}-g).
\end{aligned}
\]
Since $p-s_iq\neq 0$, we have $\theta^{-n}(\theta^k-1)(\theta^{n-m+k}-g)=0$. 
However, this equation holds only for finitely many $p$, which is a contradiction.
\end{proof}

Proposition \ref{prop:no_indefinite} leads a strong result: An indefinite component except for $Z_{\mathrm{steady}}$ must be infinite or vertical.

\subsubsection{Classification of irreducible components}

Summarizing, irreducible components of $S$ can be classified into the following four types (see Figure \ref{fig:1}):
\begin{itemize}
\item[(X)] a component that does not intersect $\{p=0\}$,
\item[($\infty$)] a component that intersects $\{p=0\}$ at $P_k$ for some $1\leq k\leq m-1$,
\item[(Steady)] the component of the steady state $Z_{\mathrm{steady}}$ (only when $g=1$),
\item[(E)] a component that intersects $\{p=0\}$ at a point in $\{a_m\neq 0\}\setminus \{P_{\mathrm{steady}}\}$.
\end{itemize}

Each component of type X is contained in a hyperplane $\{p=c^\ast\}$ for some $c^\ast\in \CC\setminus\{0\}$; hence, it is a vertical component.
We can safely ignore them when we are dealing with a generic value $p$.
(An example of type X component will be given in Appendix \S \ref{sec:app_sol_list} and \S \ref{sec:more_examples}.)

As mentioned at the beginning of \S \ref{sec:aroung_p_0}, no admissible Bethe solutions can be obtained from irreducible components of type $\infty$.

As the sole exception, $Z_{\mathrm{steady}}$ yields a Bethe vector although it is an indefinite component.
By convention, we define the Bethe vector for the steady state solution $z_1=\dots=z_m=1$ to be $[1,1,\dots,1]^T$.

We expect all admissible Bethe roots to be derived from irreducible components of type E and $Z_{\mathrm{steady}}$. 
If an irreducible components $Z$ is (i) of Type E or  (ii) $g=1$ and $Z=Z_{\mathrm{steady}}$, we will refer $Z$ to as \textit{Type ES component}.

\begin{figure}[htpb]
\centering
\begin{tikzpicture}[yscale=.8, xscale=.75]
\def\xkiten{-5}
\def\ykiten{1}
\def\yoko{10}
\def\tate{7.5}
\def\xkatamuki{2.3}
\def\ykatamuki{2}
\draw[->](-6,0)--(8,0)node[below]{$\CC$};
\draw (\xkiten,\ykiten)--++(\yoko,0)--++(0,\tate)--++(-\yoko,0)--cycle;
\draw (\xkiten,{\ykiten+\tate})--++(\yoko,0)--++(\xkatamuki,\ykatamuki)--++(-\yoko,0)node[left = 10]{$\PP^m\times \PP^{n-1}\times \CC$}--cycle;
\draw ({\xkiten+\yoko+\xkatamuki},{\ykiten+\tate+\ykatamuki})--++(0,-\tate)--
++(-\xkatamuki,-\ykatamuki);
\draw[dashed] ({\xkiten+\xkatamuki},{\ykiten+\tate+\ykatamuki})--++(0,-\tate)--
++(-\xkatamuki,-\ykatamuki);
\draw[dashed] ({\xkiten+\xkatamuki},{\ykiten+\ykatamuki})--++(\yoko,0);
\draw[->] (7.5,2)--(7.5,.8)node[above right]{$\pi$};
\draw[fill](1,0) circle[radius = 0.1]node[below]{$p=0$};
\draw[dashed] (1,0)--(-1,1)--++(\xkatamuki,\ykatamuki)--cycle;
\draw[fill, gray, fill opacity=0.2] (-1,1)--++(0,\tate)--++(\xkatamuki,\ykatamuki)--++(0,-\tate)--cycle;
\draw (-1,1)--++(0,\tate)--++(\xkatamuki,\ykatamuki)--++(0,-\tate)--cycle;
\fill[black, fill opacity=.1] (\xkiten,6)node[left, fill opacity=1]{$\{a_0=0\}$}--++(\yoko,0)--++(\xkatamuki,\ykatamuki)--++(-\yoko,0)--cycle;
\draw (\xkiten,6)--++(\yoko,0)--++(\xkatamuki,\ykatamuki);
\draw[dashed] (\xkiten,6)--++(\xkatamuki,\ykatamuki)--++(\yoko,0);
\draw[dotted,thick] (-1,6)--++(\xkatamuki,\ykatamuki);
\fill ($(-1,6)!0.3!({-1+\xkatamuki},{6+\ykatamuki})$) circle[radius = .1] node[below]{$P_1$};
\fill ($(-1,6)!0.5!({-1+\xkatamuki},{6+\ykatamuki})$) circle[radius = .1]node[above left]{$P_2$};
\fill ($(-1,6)!0.8!({-1+\xkatamuki},{6+\ykatamuki})$) circle[radius = .1] node[right = 28, above = 6]{$P_{m-1}$};
\fill (-2,7.8) circle[radius = .1];
\fill (-0.4,1.8) circle[radius = .1];
\fill (-0.1,2.2) circle[radius = .1];
\fill (0.2,3.84) circle[radius = .1];
\fill (0.8,4.66) circle[radius = .1];
\fill (-3.4,6.3) circle[radius = .1];
\fill (0.5,5) circle[radius = .1]node[above]{$P_{\mathrm{steady}}$};
\begin{scope}
\path[clip] (-8,0)--(-8,6)--(-2,6)--(-2,0)--cycle;
\draw [thick] plot [smooth, tension=2] coordinates { (-6,2.1) (4,3) (-1,4.3) (8,4.6)}node[below]{type E};
\end{scope}
\begin{scope}
\path[clip] (10,0)--(10,6)--(-1.5,6)--(-1.5,0)--cycle;
\draw [thick] plot [smooth, tension=2] coordinates { (-6,2.1) (4,3) (-1,4.3) (8,4.6)}node[below]{type E};
\end{scope}
\begin{scope}
\path[clip] (-10,4)--(-10,6)--(-3.4,6)--(-3.4,4)--cycle;
\draw[thick] (-5.5,5)node[left]{$\{\mathcal{E}=0\}\supset Z_{\mathrm{steady}}$}node[below =10, left]{(if $g=1$)}--(8,5);
\end{scope}
\begin{scope}
\path[clip] (8,4)--(8,6)--(-3.1,6)--(-3.1,4)--cycle;
\draw[thick] (-5.5,5)node[left]{$\{\mathcal{E}=0\}\supset Z_{\mathrm{steady}}$}node[below =10, left]{(if $g=1$)}--(8,5);
\end{scope}
\begin{scope}
\path[clip] (-6,7.6)--(-6,11)--(0,11)--(0,7.6)--cycle;
\draw[thick] (-2,2.6)--++(0,\tate)node[right = 6]{Type X};
\end{scope}
\begin{scope}
\path[clip] (-6,7.4)--(-6,5.1)--(0,5.1)--(0,7.4)--cycle;
\draw[thick] (-2,2.6)--++(0,\tate)node[right = 6]{Type X};
\end{scope}
\begin{scope}
\path[clip] (-6,0)--(-6,4.9)--(0,4.9)--(0,0)--cycle;
\draw[thick] (-2,2.6)--++(0,\tate)node[right = 6]{Type X};
\end{scope}
\begin{scope}
\path[clip] (-6,0)--(-6,8)--(-3.5,8)--(-3.5,0)--cycle;
\draw [thick] plot [smooth, tension=.9] coordinates { (-4,7) (-2,7.5) ($(-1,6)!0.3!({-1+\xkatamuki},{6+\ykatamuki})$) (7.8,7)}node[right]{Type $\infty$};
\end{scope}
\begin{scope}
\path[clip] (-3.4,0)--(-3.4,8)--(11,8)--(11,0)--cycle;
\draw [thick] plot [smooth, tension=.9] coordinates { (-4,7) (-2,7.5) ($(-1,6)!0.3!({-1+\xkatamuki},{6+\ykatamuki})$) (7.8,7)}node[right]{Type $\infty$};
\end{scope}
\draw [thick] plot [smooth, tension=2] coordinates { (1.4,7) ($(-1,6)!0.5!({-1+\xkatamuki},{6+\ykatamuki})$) (6,7)};
\draw [thick] plot [smooth, tension=2] coordinates { (0,8) ($(-1,6)!0.8!({-1+\xkatamuki},{6+\ykatamuki})$) (6,7.7)};
\begin{scope}
\draw [thick] plot [smooth, tension=.6] coordinates { (-3.5,8.2)(-1.9,2.2)(6.7,2)};
\end{scope}
\end{tikzpicture}
\caption{The subscheme $S\subset \PP^{m}\times \PP^{n-1}\times \CC$ is generically $1$ dimensional.
Only the irreducible components of type E can yield an admissible Bethe roots.
They may intersect the hyperplanes $\{a_0=0\}$ and $\{\mathcal{E}=0\}$ at isolated points.
}
\label{fig:1}
\end{figure}

\subsection{Degree counting}

We count points in $S|_p$ to obtain the number of quasi-simple solutions to the ASEP-Bethe equations.
Thanks to defining the problem on $\PP^m\times \PP^{n-1}$, we can solve it using only standard intersection theory.
(A brief preliminaries on intersection theory is given in Appendix \ref{sec:app:intersection}.
For readers who are interested in the intersection theory, we recommend the textbook by Fulton~\cite{fulton2013intersection}.)
Throughout the following discussion, we preferentially use terminology standard in algebraic geometry.

Let $\mathbb{A}^1_{\CC}$ be the algebraic line, whose underlying set is $\CC$, and $\pi|_S:S\to \mathbb{A}^1_\CC$ be the morphism of schemes induced from the projection \eqref{eq:projection_def}.
Define the scheme
\[
\pi_\mathcal{F}:\mathcal{F}\to \mathbb{A}^1_{\CC}
\]
over $\mathbb{A}^1_\CC$ obtained from $\pi|_S:S\to \mathbb{A}^1_\CC$ by taking the generic fiber of $\pi|S$ and its scheme-theoretic closure.
In general, $\mathcal{F}$ need not coincide with $S$. 
(In fact, any vertical component in $S$ is removed from $\mathcal{F}$.)
However, 
there exists a nonempty Zariski-open subset $U\subset \mathbb{A}^1_\CC$ such that $\mathcal{F}|_U:=\pi^{-1}_\mathcal{F}(U)$ is isomorphic to $S|_U:=(\pi|_S)^{-1}(U)$.
By construction, the morphism $\pi_\mathcal{F}$ is finite and flat (see, for example \cite[III. \S 9]{hartshorne2013algebraic}).
Note that this construction does not grantee $\mathcal{F}$ is reduced.

\begin{rem}\label{rem:0_is_in_U}
We can assume the Zariski-open set $U\subset \mathbb{A}^1_\CC$ defined in the above paragraph to contain the closed point $0$.
In fact, because $S\subset \PP^m\times \PP^{n-1}\times \mathbb{A}^1_\CC$ is an intersection of $(m+n-1)$ effective Cartier divisors and $\dim(S|_0)=0$ (Theorem \ref{thm:0dim}), $S$ is $1$-dimensional locally complete intersection around each point in $S|_0$.
Since (i) $\mathbb{A}^1_\CC$ is non-singular at $0$, (ii) $S$ is locally complete intersection (therefore, the local ring $\mathcal{O}_{S,x}$ is Cohen-Macaulay for $x\in S|_0$), and (iii) $\dim(\mathcal{O}_{S,x})=\dim(\mathcal{O}_{\mathbb{A}^1_\CC,0})+\dim(\mathcal{O}_{S|_0,x})$, there exists some open set $S|_0\subset \exists W\subset S$ such that the restriction $\pi|_{W}:W\to \mathbb{A}^1_\CC$ is flat by the ``miracle flatness theorem'' (see \cite[Theorem 23.1]{matsumura1989commutative}).
Let $V\subset \mathbb{A}^1_\CC$ be the maximal subset satisfying 
\[
(\pi|_S)^{-1}(V)\subset W.
\]
(Precisely, define $V:=\mathbb{A}^1_\CC\setminus \pi(S\setminus W)$).
As $\pi$ is proper, $V$ is a open neighborhood of $0$. 
Since $\pi|_{S_V}:S_V\to V$ is finite and flat, $\mathcal{F}_V$ and $S_V$ are isomorphic.
\end{rem}


%


In a similar manner, we define the subschemes
\[
\pi_{\mathrm{ES}}:\mathcal{F}_{\mathrm{ES}}\to \mathbb{A}^1_\CC,\qquad 
\pi_{\infty}:\mathcal{F}_{\infty}\to \mathbb{A}^1_\CC
\]
of $\mathcal{F}$ as follows: 
Let $U\subset \mathbb{A}^1_\CC$ be the Zariski-open set defined in the previous paragraph and $\mathcal{F}_{\mathrm{ES}}|_U$ and $\mathcal{F}_{\infty}|_U$ be the subschemes of $\mathcal{F}|_U$ defined as follows.
\begin{itemize}
\item $\mathcal{F}_{\mathrm{ES}}|_U$ is the unique open and closed subscheme of $\mathcal{F}|_U$ whose underlying set is the union of all irreducible components of type ES.
\item $\mathcal{F}_{\infty}|_U$ is the unique open and closed subscheme of $\mathcal{F}|_U$ whose underlying set is the union of all irreducible components of type $\infty$.
\end{itemize}
We can assume $\mathcal{F}_{\mathrm{ES}}|_U$ and $\mathcal{F}_{\infty}|_U$ are separated by removing finitely many points from $U$ if needed.
(Since these components are separated at the hyperplane $\{p=0\}$, we can retain the condition $0\in U$.)
We define $\mathcal{F}_{\mathrm{ES}}$ and $\mathcal{F}_{\infty}$ by taking the scheme-theoretic closure of $\mathcal{F}_{\mathrm{ES}}|_U$ and $\mathcal{F}_{\infty}|_U$, respectively.
$\pi_{\mathrm{ES}}$ and $\pi_\infty$ are also finite and flat.

For generic $p$, the fiber $\mathcal{F}_{\mathrm{ES}}|_p$ consists of all points in $S|_p$ contained in the components of type ES.
To count the solutions to the ASEP-Bethe equations, we focus on the number of points in $\mathcal{F}_{\mathrm{ES}}|_p$.


For an algebraic variety $X$, let $A_{k}(X)$ be the group of rational equivalence classes of $k$-cycles on $X$ (\S \ref{sec:cycles}).
Write $A^{k}(X):=A_{\dim(X)-k}(X)$. 
Let $A^\ast(X):=\bigoplus_{k=0}^{\dim(X)} A^k(X)$ be the \textit{Chow ring} of $X$ (\S \ref{sec:chow_ring}).
$A^\ast(X)$ is a graded ring with the intersection product.

For the ambient space $Y:=\PP^m\times \PP^{n-1}\times \mathbb{A}^1_\CC$,
%
we write
\[
Y_p:=\pi^{-1}(p).
\]
For a closed point $p\in \mathbb{A}^1_\CC$, we let $p:\{p\}\to \mathbb{A}^1_\CC$ denote the embedding.
Then, there exists the \textit{refined Gysin map} \cite[\S 10]{fulton2013intersection}
\[
p^!:A_{k+1}(Y)\to A_{k}(Y_p),
\]
which sends a $(k+1)$-cycle $\alpha\in A_{k+1}(Y)$ to $\alpha_p:=p!(\alpha)\in A_k(Y_p)$.
The cycle $\alpha_p$ is referred to as the \textit{specialization of $\alpha$ to $p$}.
If a cycle $\alpha=[Z]$ is defined by a subscheme $Z\subset Y$ of pure-dimension $k+1$, which admits a flat morphism $\pi|_Z:Z\to \mathbb{A}^1_\CC$,
the specialization $\alpha_p$ satisfies
\[
\alpha_p=[Z|_p]\qquad \in A_k(Y_p).
\]

Let 
\[
\eta:=[\mathcal{F}],\quad 
\eta_{\mathrm{ES}}:=[\mathcal{F}_{\mathrm{ES}}],
\quad 
\eta_{\infty}:=[\mathcal{F}_{\infty}]
\qquad \in A_{1}(Y)
\]
be the rational equivalent classes defined by $\mathcal{F}$, $\mathcal{F}_{\mathrm{ES}}$, and $\mathcal{F}_{\infty}$; then, $\eta=\eta_{\mathrm{ES}}+\eta_{\infty}$.
Since $\pi_{\mathrm{ES}}$ is finite, the total number of closed points in $\mathcal{F}_{\mathrm{ES}}|_p$ counted with multiplicity is expressed as
\[
\mathcal{N}_p:=\deg({\eta_{\mathrm{ES}}}_p),
\]
where $\deg$ is the degree of a $0$-cycle \cite[Definition 1.4]{fulton2013intersection}.

As $\pi_{\mathrm{ES}}$ is flat, we can apply the \textit{principle of continuity} \cite[Proposition 10.2]{fulton2013intersection} to conclude that $\mathcal{N}_p$ is independent of $p$.
In particular, it holds that
\[
\mathcal{N}_p=\mathcal{N}_0
\]
for all closed points $p\in \mathbb{A}^1_\CC$.

\begin{rem}
For generic $p$, $\mathcal{N}_p$ also equals the number of points in $S|_p$ contained in irreducible components of type ES counted with multiplicity.
If one desires to show that ``the number of such points is exactly $\mathcal{N}_p$ for generic $p$,'' one must assume $\mathcal{F}_{\mathrm{ES}}$ to be reduced. 
If there exists a specific value $c\in \CC$ such that $S|_c$ contains exactly $\mathcal{N}_c$ such points, the claim then true.
\end{rem}

Let $H^{a}_k=\{a_{k}=0\}$, $H^{b}_k=\{b_k=0\}$ be the hyperplanes in $\PP^m\times \PP^{n-1}$.
Denote the $(m+n-2)$-cycles defined by these hyperplanes by $[H^a_k]$ and $[H^b_k]$, respectively.
In the Chow ring $A^\ast(Y_p)=A^\ast(\PP^{m}\times \PP^{n-1})$, let $a$ and $b$ denote the common rational equivalence classes of $[H^a_0],\dots,[H^a_m]$ and $[H^b_0],\dots,[H^b_{n-1}]$, respectively.
Then, $A^\ast(Y_p)$ has the following presentation \cite[Example 8.4.2]{fulton2013intersection}
\begin{equation}\label{eq:presentation_of_AY}
A^\ast(Y_p)=\ZZ[a,b]/(a^{m+1},b^{n}).
\end{equation}
In particular, the class of a closed point $[\mathrm{pt}]$ is given by
\[
[\mathrm{pt}]=a^{m}b^{n-1}.
\]
We define $n$-cycles $\mathcal{H}_0,\mathcal{H}_1,\dots,\mathcal{H}_{m-1}\in A_n(\PP^m\times \PP^{n-1})$ by
\[
\mathcal{H}_0:=\bigcap_{i=1}^{m-1} H^b_{n-i},\quad
\mathcal{H}_k:=\bigcap_{i=1}^{m-k-1} H^b_{n-i}\cap \bigcap_{j=1}^{k} H^a_{m+1-j},
\]
whose rational equivalence classes in $\PP^{m}\times \PP^{n-1}$ are given by
\[
\mathcal{H}_k=a^kb^{m-k-1}\quad \in A_n(\PP^{m}\times \PP^{n-1}).
\]

Let $D_k$ be the Cartier divisor on $Y$ defined by $f_k$ (Eq.~\eqref{eq:def_of_f_k}), and $[D_k]\in A^1(Y)$ its rational equivalent class.
Since each $f_k$ is a homogeneous polynomials in $a_i$ and $b_i$ of bi-degree $(1,1)$, the restriction $[D_k]_p$ is expressed as
\begin{equation}\label{eq:D_is_a_add_b}
[D_k]_p=a+b\quad \in A^{1}(Y_p)
\end{equation}
for all $p\in \mathbb{A}^1_\CC$.

Let $U\subset \mathbb{A}^1_\CC$ be an Zariski-open set such that $\mathcal{F}|_U$ and $S|_U$ are isomorphic.
We can assume $0\in U$ (see Remark \ref{rem:0_is_in_U}).
Then, we have
\begin{equation}\label{eq:eta_0_expression}
\eta_0=[S|_0]=[D_1]_0\cdot [D_2]_0\cdot \cdots \cdot [D_{m+n-1}]_0\qquad \in A_{0}\left(S|_0\right).
\end{equation}


We want to calculate the subclass $(\eta_{\mathrm{ES}})_0$ of $\eta_0$.
It is straightforward to obtain
\begin{equation}\label{eq:inetersection_n_to_m-1}
[D_{n+1}]_0\cdot [D_{n+2}]_0\cdot \dots\cdot [D_{m+n-1}]_0=(a+b)^{m-1}
=\sum_{k=0}^{m-1} 
\binom{m-1}{k}
 a^kb^{m-k-1}
\quad \in A^{m-1}(Y_0).
\end{equation}
Since $f_k(a,b;0)$ for $k\geq n+1$ has the simple expression
\[
f_k(a,b;0)=\sum_{
\substack{i+j=k
\\i\geq 0,j\geq 1
}} 
a_ib_j,
\]
the intersection of the supports $|D_{n+1}|_0|,\dots,|D_{m+n-1}|_0|$ is expressed as
\begin{equation}\label{eq:hyper_inf_planes}
|D_{n+1}|_0|\cap |D_{n+2}|_0|\cap \dots\cap |D_{m+n-1}|_0|
= \bigcup_{k=0}^{m-1}|\mathcal{H}_k|.
\end{equation}
As the correspondence $\mathcal{H}_k\mapsto a^kb^{m-k-1}$ ($k=0,1,\dots,m-1$) is one-to-one, and $a^kb^{m-k-1}$ ($k=0,1,\dots,m-1$) are linearly independent in $A_n(Y_0)$, the equation \eqref{eq:inetersection_n_to_m-1} can be uniquely lifted to
the group $A_{n}(\textstyle \bigcup_{k=0}^{m-1}|\mathcal{H}_k|)$ as
\begin{equation}\label{eq:cycle_D_n}
[D_{n+1}]_0\cdot [D_{n+2}]_0\cdot \dots\cdot [D_{m+n-1}]_0
=\sum_{k=0}^{m-1}\binom{m-1}{k}\cdot \mathcal{H}_k\quad \in A_{n}(\textstyle \bigcup_{k=0}^{m-1}|\mathcal{H}_k|). 
\end{equation}

Next consider the class
\begin{equation}\label{eq:inter_of_H}
[D_{1}]_0\cdot [D_{2}]_0\cdot\cdots \cdot [D_{n}]_0\cdot \mathcal{H}_k
=(a+b)^na^kb^{m-k-1}=\binom{n}{m-k}\cdot [\mathrm{pt}]\quad \in A^{m+n-1}(Y_0).
\end{equation}
If $k>0$, we have the set-theoretical equation
\begin{equation}\label{eq:inetersection_D_1_to_n}
|D_{1}|_0|\cap |D_{2}|_0|\cap \dots\cap |D_{n}|_0| \cap  |\mathcal{H}_k|
= \{P_k\},
\end{equation}
which follows from Proposition \ref{prop:bound_points} and \eqref{eq:hyper_inf_planes}.
Then, it holds from \eqref{eq:inter_of_H} that
\begin{equation}\label{eq:aux_2}
\begin{aligned}
&[D_{1}]_0\cdot [D_{2}]_0\cdot\cdots \cdot [D_{n}]_0\cdot \mathcal{H}_k=\binom{n}{m-k}\cdot P_k\quad &\in A^{0}(P_k).
\end{aligned}
\end{equation}

Let
\[
K=|D_1|_0|\cap |D_2|_0|\cap \dots\cap |D_{m+n-1}|_0|
=\bigcup_{k=0}^{m-1} |D_1|_0|\cap |D_2|_0|\cap \dots\cap |D_{n}|_0|\cap |\mathcal{H}_k|.
\]
By \eqref{eq:inetersection_D_1_to_n} and $\mathcal{H}_0\cap \{P_1,\dots,P_{m-1}\}=\emptyset$, $A^\ast(K)$ is decomposed as 
\begin{equation}\label{eq:direct_decomp}
A^\ast(K)=A^\ast(|D_1|_0|\cap |D_2|_0|\cap \dots\cap |D_{n}|_0|\cap |\mathcal{H}_0|)\oplus \bigoplus_{k=1}^{m-1} A^\ast(P_k).
\end{equation}
In $A^\ast(K)$, it follows from \eqref{eq:eta_0_expression}, \eqref{eq:cycle_D_n} and \eqref{eq:aux_2} that
\[
\begin{aligned}
\eta_0=[D_1]_0\cdot [D_2]_0\cdots [D_{m+n-1}]_0
&=
\sum_{k=0}^{m-1}\binom{m-1}{k}[D_1]_0\cdot [D_2]_0\cdots [D_{n}]_0\cdot\mathcal{H}_k\\
&=
[D_1]_0\cdot [D_2]_0\cdots [D_{n}]_0\cdot\mathcal{H}_0+
\sum_{k=1}^{m-1}\binom{m-1}{k}\binom{n}{m-k}\cdot P_k.
\end{aligned}
\]
Comparing this equation and the direct decomposition \eqref{eq:direct_decomp}, we obtain
\[
(\eta_{\mathrm{ES}})_0=
[D_1]_0\cdot [D_2]_0\cdots [D_{n}]_0\cdot\mathcal{H}_0,\qquad
(\eta_\infty)_0=\sum_{k=1}^{m-1}\binom{m-1}{k}\binom{n}{m-k}\cdot P_k.
\]
Therefore, it follows from \eqref{eq:inter_of_H} that
\[
\mathcal{N}_p=
\mathcal{N}_0=\deg((\eta_{\mathrm{ES}})_0)=
\deg\left([D_1]_0\cdot [D_2]_0\cdots [D_{n}]_0\cdot\mathcal{H}_0\right)=
\binom{n}{m}.
\]

Summering, we obtain the following theorem:
\begin{theorem}
For generic $p$, the number of points in $S|_p$ contained in irreducible components of type ES is $\binom{n}{m}$ with multiplicity.
If $\mathcal{F}_{\mathrm{ES}}$ is reduced, the number of such points is exactly $\binom{n}{m}$.
\end{theorem}

\begin{cor}\label{cor:conclusion}
For generic $g$ and $p$, the number of points in $S|_p$ contained in irreducible components of type ES is exactly $\binom{n}{m}$.
\end{cor}
\begin{proof}
It suffices to show that the number of points in $S|_0$ in irreducible components of type ES is exactly $\binom{n}{m}$ for generic $g$.
This fact is proved by Iwao-Motegi~\cite[Corollary 2.8]{iwao2025bethe}.
\end{proof}

By Corollary \ref{cor:conclusion}, we obtain the following conclusion:
For generic $g$ and $p$, the number of quasi-simple solutions to the ASEP-Bethe equations is $\binom{n}{m}$.

\subsection{Examples}

We give small examples to describe our method.
Throughout, we assume $g\neq 0$.

\subsubsection{$m=2$ case}

Suppose $m=2$.
Let $\{j,k\}=\{1,2\}$.
In this case, we have
\[
\begin{aligned}
\varphi(z_j;p)&=z_j^n(qz_kz_j-z_k+p)\\
&=z_j^n(qe_2-z_j^{-1}e_2+p(z_j^{-1}e_1-z_j^{-2}e_2))
=qe_2z_j^{n}+(pe_1-e_2)z_j^{n-1}-pe_2z^{n-2}
\end{aligned}
\]
and
\[
\begin{aligned}
\psi(z_j;p)
&=qz_kz_j-z_j+p
=q(e_1-z_j)z_j-z_j+p
=-qz_j^2+(qe_1-1)z_j+p.
\end{aligned}
\]
Then, the projectivized problem is
\begin{equation}\label{eq:example_n_2}
\begin{bmatrix}
a_0b_1\\
a_0b_2+a_1b_1\\
a_0b_3+a_1b_2+a_2b_1\\
\vdots\\
\vdots\\
a_1b_{n-1}+a_2b_{n-2}\\
a_2b_{n-1}
\end{bmatrix}
=b_0
\begin{bmatrix}
qa_2\\
a_2-pa_1\\
-pa_2\\
0\\
\vdots\\
0\\0
\end{bmatrix}
-(-1)^n
b_0g
\begin{bmatrix}
0\\
0\\
0\\
\vdots\\
qa_0\\
qa_1-a_0\\
-pa_0
\end{bmatrix}.
\end{equation}

\begin{example}\label{ex:m=2_exaple}
If $n=4$, we have
\begin{equation}\label{eq:42example}
\begin{bmatrix}
a_0b_1\\
a_0b_2+a_1b_1\\
a_0b_3+a_1b_2+a_2b_1\\
a_1b_3+a_2b_2\\
a_2b_3
\end{bmatrix}
=b_0
\begin{bmatrix}
qa_2\\
a_2-pa_1\\
-pa_2\\
0\\0
\end{bmatrix}
-
b_0g
\begin{bmatrix}
0\\
0\\
qa_0\\
qa_1-a_0\\
-pa_0
\end{bmatrix}.
\end{equation}
Let $S$ be the solution set of the system \eqref{eq:42example} and $(a,b):=([a_0:a_1:a_2],[b_0:b_1:b_2:b_3])$ be the standard coordinate of points in $S$.

It is straightforward to show that $a_0=0\iff a_2=0$.
If $a_0=0$, the only solution to \eqref{eq:42example} is
\begin{equation}\label{eq:example_sol_1_1}
(a,b)=([0:1:0],[1:-p:0:-gq]).
\end{equation}
If $a_0\neq 0$, one can normalize $a_0=b_0=1$.
In this case, we have $b_1=qa_2$ and $b_2=-qa_1a_2+a_2-pa_1$ from the first and the second equation in the system \eqref{eq:42example}.
Calculating 
\[
a_2\times (\text{the fourth equation})-a_1\times (\text{the fifth equation})
\] 
yields $(a_2^2-g)b_2=0$.
Using this, we soon derive the following six solutions for a generic $p$:
\begin{equation}\label{eq:example_sol_1_2}
\begin{aligned}
(a,b)
=&([1:\frac{\sqrt{-g}}{q\sqrt{-g}\pm p}:\pm\sqrt{-g}],[1:\pm q\sqrt{-g}:0:\mp p\sqrt{-g}]),\\
&([1:\frac{\sqrt{g}\pm \sqrt{g+8\sqrt{g}D^2_+}}{2D_+}:\sqrt{g}],[1:q\sqrt{g}:\frac{\sqrt{g}\mp \sqrt{g+8\sqrt{g}D^2_+}}{2}:p\sqrt{g}]),\\
&([1:\frac{-\sqrt{g}\pm \sqrt{g-8\sqrt{g}D^2_-}}{2D_-}:-\sqrt{g}],[1:-q\sqrt{g}:\frac{-\sqrt{g}\mp \sqrt{g-8\sqrt{g}D^2_-}}{2}:-p\sqrt{g}]),
\end{aligned}
\end{equation}
where $D_\pm=p\pm q\sqrt{g}$.
The family \eqref{eq:example_sol_1_1} forms an irreducible component of type $\infty$, and families in \eqref{eq:example_sol_1_2} form irreducible components of type ES.
Note that the components of type E may intersect $\{a_0=0\}$ at an isolated point.
\end{example}

Note that the ``$b$-coordinates'' $[b_0:b_1:\dots:b_{n-1}]$ of solutions to \eqref{eq:example_n_2} are uniquely determined from their $a$-coordinates by
\[
\begin{bmatrix}
b_1\\b_2\\\vdots\\\vdots \\b_{n-1}
\end{bmatrix}
=
\begin{bmatrix}
a_0\\
a_1&a_0\\
a_2&a_1&a_0\\
&\ddots&\ddots&\ddots\\
&&a_2&a_1&a_0
\end{bmatrix}^{-1}
\cdot 
\left(b_0
\begin{bmatrix}
qa_2\\
a_2-pa_1\\
-pa_2\\
0\\\vdots\\0
\end{bmatrix}
-
b_0g
\begin{bmatrix}
0\\\vdots\\
0\\
qa_0\\
qa_1-a_0\\
-pa_0
\end{bmatrix}
\right)
\]
whenever $a_0\neq 0$.
Therefore, for brevity, we will often omit the $b$-coordinates in the sequel.

From \eqref{eq:example_sol_1_2}, we obtain a complete list of quasi-simple solutions to ASEP-Bethe equations for $(m,n)=(2,4)$.
In particular, for the case $g=1$, we obtain
\begin{equation}\label{eq:six_solutions}
\begin{aligned}
[a_0:a_1:a_2]=
&[1:\frac{\sqrt{-1}}{q\sqrt{-1}\pm p}:\pm\sqrt{-1}],
[1:2:1],[1:-1:1],
[1:\frac{-1\pm \sqrt{1-8(p-q)^2}}{2(p-q)}:-1].
\end{aligned}
\end{equation}
From these, we obtain the following $6=\binom{4}{2}$ solutions to the ASEP Bethe equations \eqref{eq:ASEP_Bethe}:
\[
\{z_1,z_2\}=\{c^{+}_1,c^{+}_2\},\{c^{-}_1,c^{-}_2\},\{1,1\},\{\frac{-1+\sqrt{-3}}{2},\frac{-1-\sqrt{-3}}{2}\},\{d^+_1,d^+_2\},\{d^-_1,d^-_2\},
\]
where $c^\pm_1,c^{\pm}_2$ are roots of $c^2-\frac{\sqrt{-1}}{q\sqrt{-1}\pm p}c\pm \sqrt{-1}$, and 
$d^\pm_1,d^{\pm}_2$ are roots of $d^2-\frac{-1\pm \sqrt{1-8(p-q)^2}}{2(p-q)}d-1$.
It is straightforward to verify that they satisfy \eqref{eq:ASEP_Bethe}, except for the steady state solution $\{1,1\}$.

\subsubsection{$m=3$ case}

Suppose $m=3$. 
Let $\{j,k,l\}=\{1,2,3\}$.
We let $e_i$ (\textit{resp}. $e_i'$) denote the $i$-th elementary symmetric polynomial in $z_1,z_2,z_3$ (\textit{resp}. $z_k,z_l$).
Then, we have
\[
\begin{aligned}
\varphi(z_j;p)
&=z_j^n(qz_kz_j-z_k+p)(qz_lz_j-z_l+p)\\
&=z_j^n\{(qz_j-1)^2e_2'+p(qz_j-1)e_1'+p^2\}\\
&=q^2 e_3 z_j^{n+1} + (p q e_2 - 2 q e_3) z_j^n + \bigl(p^2 e_1 - p e_2 + (1 - p q) e_3\bigr) z_j^{n-1} \\
&\hspace{16em} + (p e_3 - p^2 e_2) z_j^{n-2} + p^2 e_3 z_j^{n-3}
\end{aligned}
\]
and
\[
\begin{aligned}
\psi(z_j;p)
&=(qz_kz_j-z_j+p)(qz_lz_j-z_j+p)\\
&=(qz_k-1)(qz_l-1)z_j^2+\{(qz_k-1)+(qz_l-1)\}z_jp+p^2\\
&=(q^2e_2'-qe_1'+1)z_j^2+(qe_1'-2)z_jp+p^2\\
&=q^2 z_j^4 + q(1 - q e_1) z_j^3 + (q^2 e_2 - q e_1 - p q + 1) z_j^2 + p(q e_1-2) z_j + p^2.
\end{aligned}
\]

Therefore, the projectivized system is expressed as follows:
\begin{equation}\label{eq:n3example}
\begin{bmatrix}
a_0b_1\\
a_0b_2+a_1b_1\\
a_0b_3+a_1b_2+a_2b_1\\
\vdots\\
\vdots\\
a_2b_{n-1}+a_3b_{n-2}\\
a_3b_{n-1}
\end{bmatrix}
=b_0
\begin{bmatrix}
q^2a_3\\
2 q a_3-p q a_2\\
p^2 a_1 - p a_2 + (1 - p q) a_3\\
p^2 a_2-p a_3\\
p^2 a_3\\
0\\\vdots\\0
\end{bmatrix}
+(-1)^n
b_0g
\begin{bmatrix}
0\\\vdots\\0\\
q^2a_0\\
q^2a_1-qa_0\\
q^2a_2-qa_1+(1-pq)a_0\\
2pa_0-pqa_1\\
p^2a_0
\end{bmatrix}.
\end{equation}

\begin{example}\label{ex:example_4_3}
If $n=4$, we have
\begin{equation}\label{eq:43example}
\begin{bmatrix}
a_0b_1\\
a_0b_2+a_1b_1\\
a_0b_3+a_1b_2+a_2b_1\\
a_1b_3+a_2b_2+a_3b_1\\
a_2b_3+a_3b_2\\
a_3b_3
\end{bmatrix}
=b_0
\begin{bmatrix}
q^2a_3\\
2 q a_3-p q a_2\\
p^2 a_1 - p a_2 + (1 - p q) a_3\\
p^2 a_2-p a_3\\
p^2 a_3\\
0
\end{bmatrix}
+
b_0g
\begin{bmatrix}
0\\
q^2a_0\\
q^2a_1-qa_0\\
q^2a_2-qa_1+(1-pq)a_0\\
2pa_0-pqa_1\\
p^2a_0
\end{bmatrix}.
\end{equation}

Although solving the system \eqref{eq:43example} completely is not easy, one can find some special solutions by direct calculations.
(A complete list of solutions delivered by computer-assisted calculation will be given in \S \ref{sec:app_sol_list}.)

Let $S$ be the solution set of the system \eqref{eq:43example} and $(a,b):=([a_0:a_1:a_2:a_3],[b_0:b_1:b_2:b_3])$ be the standard coordinate of points in $S$.
Note that $a_0=0\iff a_3=0$.

We can find an irreducible component of $S$ of type $\infty$ by putting $a_0=0$.
In fact, if $a_0=0$, the $a$-coordinate can be parameterized as
\[
[a_0:a_1:a_2:a_3]=[0:s:t:0],\qquad [s:t]\in \PP^1.
\]
Substituting it to \eqref{eq:43example} derives
\[
pt^3(qt-s)+gqs^3(t-ps)=0
\]
and
\begin{equation}\label{eq:example_2_1}
\begin{aligned}
&[b_0:b_1:b_2:b_3]\\
&=\begin{cases}
[1:-gpq(s/t)^3+gq(s/t)^2-p:
gpq(s/t)^2-gq(s/t)+p^2+gq^2:-gpq(s/t)
] & (t\neq 0),\\
[1:-pq(t/s):pq(t/s)^2-p(t/s)+p^2+gq^2:-pq(t/s)^3+p(t/s)^2-gq] & (s\neq 0).
\end{cases}
\end{aligned}
\end{equation}
The family \eqref{eq:example_2_1} forms an irreducible component of type $\infty$.
In particular if $p=0$, we have (i) $(a,b)=([0:0:1:0],[1:0:g:0])$ or (ii) $(a,b)=([0:1:0:0],[1:0:g:-g])$.
The first solution associated to $P_1$ and the second to $P_2$.

%
\end{example}

\section{Reduction to $p=q=\frac{1}{2}$ and $g=1$}\label{sec:future}


When $p=q=\frac{1}{2}$ and $g=1$, the ASEP reduces to the SEP, which is equivalent to the Heidelberg XXX chain.
Since the rational functions representing Bethe vectors may diverge, we cannot always simply substitute the value. 
In such cases, we can introduce a perturbation along $S_{\mathrm{ess}}$. 
Since $S_{\mathrm{ess}}$ is $1$-dimensional, we can unambiguously determine the direction of the perturbation whenever $S_{\mathrm{ess}}$ is non-singular.
(If $S_{\mathrm{ess}}$ is singular, we need to apply a desingularization).

As a working example, we introduce the $(m,n)=(2,4)$ case in Example \ref{ex:m=2_exaple}, where \eqref{eq:six_solutions} shows the $6=\binom{4}{2}$ solutions to the ASEP-Bethe equations.

Let $\epsilon:=p-q$ and $i=\sqrt{-1}$.
When taking the limit $p\to \frac{1}{2}$ ($\epsilon\to 0$), we obtain the six points on $S|_{\frac{1}{2}}$ by direct calculations:
\[
\begin{gathered}
[1:\frac{i}{qi\pm p}:\pm i]\to [1:\frac{2i}{i\pm 1}:\pm i],\quad
[1:2:1]\to[1:2:1],\quad
[1:-1:1]
\to
[1:-1:1],
\\
[1:\frac{-1+ \sqrt{1-8\epsilon ^2}}{2\epsilon}:-1]\to [1:0:-1],\quad
[1:\frac{-1- \sqrt{1-8\epsilon^2}}{2\epsilon}:-1]\to [0:1:0].
\end{gathered}
\]
Each point corresponds to the solution
\begin{equation}\label{eq:solutions_example_24}
\{z_1,z_2\}=\{1,i\},\{1,-i\},\{1,1\},\{\frac{-1+i\sqrt{3}}{2},\frac{-1-i\sqrt{3}}{2}\},\{1,-1\},\{0,\infty\}.
\end{equation}
It is straightforward to check that the Bethe vectors \eqref{eq:formula_eigen_vec} are well-defined and non-zero, except for the steady state solution $\{z_1,z_2\}=\{1,1\}$ and the sixth solution $\{z_1,z_2\}=\{0,\infty\}$.

Although the formula \eqref{eq:formula_eigen_vec} is ill-defined for the steady state solution, by convention we define it as the constant vector $[1,1,1,1,1,1]^T\in \CC^{\binom{4}{2}}$.

For $\{z_1,z_2\}=\{0,\infty\}$, we can derive necessary information by perturbation along $S_{\mathrm{ess}}$.
Let $P_\ast$ be the point in $S_{\mathrm{ess}}|_{\frac{1}{2}}$ that corresponds to that solution.
Then, for small $|\epsilon|=|p-q|$, the homogeneous coordinates $[a_0:a_1:a_2]$ around $P_\ast$ along $S_{\mathrm{ess}}$ are expanded as
\[
[1:\frac{-1- \sqrt{1-8\epsilon^2}}{2\epsilon}:-1]=[\epsilon:-1+2\epsilon^2+4\epsilon^4+16\epsilon^6+\cdots:-\epsilon],
\]
which leads
\[
\{z_1,z_2\}=\{\epsilon+\epsilon^3+4\epsilon^5+\cdots,
-\epsilon^{-1} +\epsilon+3\epsilon^3+\cdots
\}.
\]

By using the Laurent expansions of $z_1$ and $z_2$, we derive the corresponding eigenvalue \eqref{eq:eigenvalue} as
\[
\lambda=\sum_{i=1}^2(qz_i-1+pz_i^{-1})=\sum_{i=1}^2\left(\left(\frac{1}{2}-\frac{\epsilon}{2}\right)z_i-1+\left(\frac{1}{2}+\frac{\epsilon}{2}\right)z_i^{-1}\right)=-\frac{3-\sqrt{1-8\epsilon^2}}{2}
\stackrel{\epsilon\to 0}{\to}-1,
\]
and the Bethe vector \eqref{eq:formula_eigen_vec} as
\[
\begin{bmatrix}
\varphi(1,2)\\
\varphi(1,3)\\
\varphi(1,4)\\
\varphi(2,3)\\
\varphi(2,4)\\
\varphi(3,4)
\end{bmatrix}
=
\begin{bmatrix}
-\epsilon^{-3}+4\epsilon^{-1}+\dots\\
-2\epsilon^{-2}+4+\dots\\
\epsilon^{-3}-4\epsilon^{-1}+\dots\\
\epsilon^{-3}-4\epsilon^{-1}+\dots\\
2\epsilon^{-2}-4+\dots\\
-\epsilon^{-3}+4\epsilon^{-1}+\dots
\end{bmatrix}
\propto
\begin{bmatrix}
-1+4\epsilon^{2}+\dots\\
-2\epsilon+4\epsilon^3+\dots\\
1-4\epsilon^2+\dots\\
1-4\epsilon^2+\dots\\
2\epsilon-4\epsilon^3+\dots\\
-1+4\epsilon^{2}+\dots
\end{bmatrix}
\stackrel{\epsilon\to 0}{\to}
\begin{bmatrix}
-1\\0\\1\\1\\0\\-1
\end{bmatrix},
\]
which is in fact an eigenvector of the Markov matrix (see \eqref{eq:spec_prob})
\[
M =
\begin{bmatrix}
-1           & \tfrac{1}{2} & 0            & 0            & \tfrac{1}{2} & 0            \\
\tfrac{1}{2} & -2           & \tfrac{1}{2} & \tfrac{1}{2} & 0            & \tfrac{1}{2} \\
0            & \tfrac{1}{2} & -1           & 0            & \tfrac{1}{2} & 0            \\
0            & \tfrac{1}{2} & 0            & -1           & \tfrac{1}{2} & 0            \\
\tfrac{1}{2} & 0            & \tfrac{1}{2} & \tfrac{1}{2} & -2           & \tfrac{1}{2} \\
0            & \tfrac{1}{2} & 0            & 0            & \tfrac{1}{2} & -1
\end{bmatrix}
\]
belonging to the eigenvalue $\lambda=-1$.

Under the variable change $z=\frac{u+\frac{i}{2}}{u-\frac{i}{2}}$, the solutions \eqref{eq:solutions_example_24} to the SEP can be transformed to solutions to the XXX model:
\[
\{u_1,u_2\}=\{\infty,\frac{1}{2}\},\{\infty,-\frac{1}{2}\},\{\infty,\infty\},\{\frac{1}{\sqrt{12}},-\frac{1}{\sqrt{12}}\},\{\infty,0\},\{-\frac{i}{2},\frac{i}{2}\}.
\]
Since the Hamiltonian of XXX and the Markov matrix $M$ is similar (see \cite{de2006exact,crampe2010eigenvectors}), all eigenvalues of XXX-Hamiltonian can be obtained from our calculations.

In Figure \ref{fig:list_of_solutions}, all solutions to the ASEP-Bethe equations corresponding to points in $S_{\mathrm{ess}}|_{\frac{1}{2}}$ are listed.
Since the number of linearly independent Bethe vectors are sufficient to give a basis of $\CC^{\binom{4}{2}}$, this is in fact a complete list of solutions.
The corresponding Bethe roots for XXX model are also listed.
Solutions A,B,C, and E are XXX-Bethe roots that contain the infinite point $\infty$.

It is pointed out by previous researches (see for example \cite{deguchi2001non}) that such ``infinite solutions'' correspond to non-highest weight vectors of the XXX-Hamiltonian, while finite solutions corresponds to highest-weight vectors.

\begin{figure}[htbp]
\centering
\begin{tabular}{c||c|c|c|c|c}
 & $\{z_1,z_2\}$ & eigenvalue& Bethe vector & $\{u_1,u_2\}$ & XXX model\\\hline
A &  $\{1,1\}$ & $0$ & $[1,1,1,1,1,1]^T$& $\{\infty,\infty\}$ & contains $2$ $\infty$'s\\
B & $\{1,i\}$ & $-1$ & $[-1+i, 0, 1+i, -1-i, 0, 1-i]^T $& $\{\infty,\frac{1}{2}\}$ & contains $1$ $\infty$\\
C & $\{1,-i\}$ & $-1$ & $[-1-i, 0, 1-i, -1+i, 0, 1+i]^T $& $\{\infty,-\frac{1}{2}\}$ & contains $1$ $\infty$\\
D & $\{0,\infty\}$ & $-1$ & $\propto[-1, 0, 1, 1, 0, -1]^T $& $\{-\frac{i}{2},\frac{i}{2}\}$ &
\!\!\!\begin{tabular}{c}
finite (physical\\
singular)
\end{tabular}
\!\!\!
\\
E & $\{1,-1\}$ & $-2$ & $[0,-2,0,0,2,0]^T$ & $\{\infty,0\}$ & contains $1$ $\infty$\\
F & $\{\omega,\omega^2\}$ & $-3$ & $[-\omega,2\omega,-\omega,-\omega,2\omega,-\omega]^T$ & $\{\frac{1}{\sqrt{12}},-\frac{1}{\sqrt{12}}\}$ & finite
\end{tabular}
\caption{A list of all quasi-simple solutions to the ASEP-Bethe equations for $(m,n)=(2,4)$.
Here, $\omega$ is a primitive $3$-rd rood of $1$.
There exists $2$ finite solutions to the XXX-Bethe equations (one of them is physical singular).
}
\label{fig:list_of_solutions}
\end{figure}

\subsection*{Acknowledgment}

This work is partially supported by KAKENHI 22K03239 and 23K03056.

\subsection*{Declaration of generated AI}

During the preparation of this manuscript, the author used ChatGPT to translate and improve the English text. 
The author has reviewed the content and takes full responsibility for the accuracy and integrity of this paper.

\appendix

\section{Preliminaries on intersection theory}\label{sec:app:intersection}

In this appendix, we give a brief review on intersection theory.
Materials and notations here are followed to the standard textbook~\cite{fulton2013intersection}.
Throughout, we work on the base field $\CC$.

\subsection{Cycles and divisors}\label{sec:cycles}

Let $X$ be an algebraic scheme.
A \textit{$k$-cycle on $X$} \cite[\S 1.3]{fulton2013intersection} is a finite formal sum
\[
\sum_i n_i\cdot [V_i],\qquad (n_i\in \ZZ),
\]
where $V_i$ is a $k$-dimensional subvariety of $X$.
All $k$-cycles form an Abelian group $Z_k(X)$.

For a cycle $\alpha=\sum_i n_i\cdot [V_i]$, the \textit{support} of $\alpha$, written $|\alpha|$ or $\mathrm{supp}(\alpha)$, is the union $\bigcup_{n_i\neq 0}V_i$ of subvarieties appearing in $\alpha$ with non-zero coefficient.
Regarding the support $|\alpha|$ as a subscheme of $X$, $\alpha$ can be seen as a cycle on $|\alpha|$; $\alpha\in Z_k(|\alpha|)$.

For a $(k+1)$-dimensional subvariety $W\subset X$ and a non-zero rational function $f$ on $W$, define the $k$-cycle $[\mathrm{div}(f)]\in Z_k(X)$ by
\[
[\mathrm{div}(f)]:=\sum_{V\subset W:\mathrm{\,codim.\, one}} \ord_{V}(f)\cdot [V],
\]
where $\mathrm{ord}_V$ is the order function $R(W)^\times\to \ZZ$ defined by the local ring $\mathcal{O}_{V,W}$.

Let $\mathrm{Rat}_k(X)\subset Z_k(X)$ be the subgroup generated by all $k$-cycles of the form $[\mathrm{div}(f)]$.
We say two $k$-cycles $\alpha_1,\alpha_2$ are \textit{rationally equivalent} if $\alpha_1-\alpha_2\in \mathrm{Rat}_k(X)$.
Define 
\[
A_k(X):=Z_k(X)/\mathrm{Rat}_k(X).
\]
Let $A_\ast(X)$ denote the direct sum of $A_k(X)$ for $k=0,1,\dots,\dim (X)$.

For a proper map $f:X\to Y$ of schemes, there exists the \textit{push-forward map} \cite[\S 1.4]{fulton2013intersection}
\[
f_\ast :A_k(X)\to A_k(Y).
\]
Let $Y_1,\dots,Y_n$ be closed subschemes of a scheme $Y$, and $\iota_i:Y_i\to Y$ be the closed embedding.
If there exists a family of classes $\alpha_i\in A_\ast(Y_i)$ and $\alpha\in A_\ast(Y)$ satisfying $\alpha=\sum_{i=1}^n(\iota_i)_\ast(\alpha_i)$, we usually use the abbreviation
\[
\alpha=\sum_{i=1}^n\alpha_i\qquad \in A_\ast(Y). 
\]
In particular, if $\iota:Y\hookrightarrow X$ is a closed embedding, we often abbreviate $i^\ast(\alpha)$ as $\alpha$ for $\alpha\in A_k(Y)$.

For a Cartier divisor\footnote{A \textit{Cartier divisor on $X$} is defined by a data $D=(U_\alpha,f_\alpha)_{\alpha}$, where $\{U_\alpha\}_{\alpha}$ is an open covering of $X$ and $f_\alpha$ is a rational function on $U_\alpha$, subject to the condition that $f_\alpha/f_\beta$ is a unit on $U_\alpha\cap U_\beta$.
If $f_\alpha/f'_\alpha$ is a unit for all $\alpha$, the two datum $(U_\alpha,f_\alpha)_{\alpha}$ and $(U_\alpha,f'_\alpha)_{\alpha}$ define the same Cartier divisor.
} $D=(U_\alpha,f_\alpha)_\alpha$ on $X$, denote the \textit{associated Weil divisor} by
\[
[D]:=\sum_{V\subset X:\mathrm{\, codim.\, one}}\mathrm{ord}_{V}(f_\alpha) \cdot [V]\quad \in Z_{k-1}(X),
\]
where $\mathrm{ord}_{V}(f_\alpha)$ is defined for $\alpha$ with $V\cap U_\alpha\neq \emptyset$.
If every $f_\alpha$ is a regular function on $U_\alpha$, $D$ is called \textit{effective}.
Any effective divisor determines a subscheme, which is defined locally by $f_\alpha$.

%

For two Cartier divisors $D_1=(U_\alpha,f_\alpha)_\alpha$ and $D_2=(U_\alpha,g_\alpha)_\alpha$, define $D_1+D_2=(U_\alpha,f_\alpha g_\alpha)_\alpha$.
Then, all Cartier divisors form an Abelian group $\mathrm{Div}(X)$ and the correspondence $D\mapsto [D]$ defines the homomorphism
\[
\mathrm{Div}(X)\to Z_{k-1}(X).
\]

\subsection{Chow ring}\label{sec:chow_ring}

For a $n$-dimensional nonsingular variety $X$, denote $A^k(X):=A_{n-k}(X)$. 
Set $A^\ast(X):=\bigoplus_{k=0}^{n}A^k(X)$.

Let $\alpha_1,\dots,\alpha_r$ be cycles on $X$ with supports $|\alpha_1|,\dots,|\alpha_r|$.
We let 
\[
\alpha_1\cdot \cdots\cdot \alpha_r \in A^\ast(|\alpha_1|\cap \dots \cap |\alpha_r|)
\]
denote the \textit{intersection class} of $\alpha_1,\dots,\alpha_r$.

There exists a bilinear homomorphism
\[
A^p(X)\otimes A^q(X)\to A^{p+q}(X)
\]
called the \textit{global intersection product} \cite[\S 8.3]{fulton2013intersection}, where the image of $\alpha\otimes \beta$ is the image of the intersection class $\alpha\cdot \beta$ under the embedding $|\alpha|\cap |\beta|\to X$.
With this product, $A^\ast(X)$ is a commutative graded ring, called as the \textit{Chow ring} of $X$.

%

\begin{example}
Let $X=\PP^M\times \PP^N$ and $([a_0:\dots:a_M],[b_0:\dots:b_N])$ be the homogeneous coordinate. 
Denote the hyperplanes $\{a_k=0\}$ and $\{b_k=0\}$ by $H^a_k$ and $H^b_k$, respectively.
As an element of $A^1(X)$, we have $a:=[H^a_0]=\dots=[H^a_M]$ and $b:=[H^b_0]=\dots=[H^b_N]$.
The Chaw ring of $X$ satisfies $A^\ast(X)=\ZZ[a,b]/(a^{M+1},b^{N+1})$.
In particular, a class of a closed point $[\mathrm{pt}]$ is given by $[\mathrm{pt}]=a^Mb^N$.
\end{example}

A scheme $Y$ is called \textit{pure dimensional} if its irreducible components has the same dimensions.

\section{A complete list of solutions to the ASEP-Bethe equations for $(m,n)=(3,4)$}\label{sec:app_sol_list}

In this section, we give a complete list of solutions to the system \eqref{eq:43example} delivered by computer-assisted calculations.

\subsubsection*{Type ES components}

Let $k_j$ $(j=0,1,2,3)$ be a arbitrary solution of $k^4-g=0$.
Define
$
\Delta_j:=(p+qk_j^2)(p^2+q^2g)+pqk_j^3
$.
Then, the irreducible components of type E are formed by the following solutions:
\[
\begin{aligned}
[a_0:a_1:a_2:a_3]
&=[1:k_j+\frac{k_j^2(p^2+pk_j+pqk_j^2+2q^2g)}{\Delta_j}:
k_j^2+\frac{k_j^3(2p^2+pqk_j^2+qk_j^3+q^2g)}{\Delta_j}:k_j^3
].
\end{aligned}
\]
with 
\[
\Delta(k)=(p+qk^2)(p^2+q^2k^4)+pqk^3.
\]

\subsubsection*{Type X components}

There are two families of solutions belonging to type X components:
\[
\begin{aligned}
[a_0:a_1:a_2:a_3]=&[1:t:k_j(k_j-1)^2-k_jt:k_j^3],\qquad p=\frac{k_j^2}{1+ k_j^2},\quad (t\in \CC),\qquad (j=0,1,2,3)\\
&[1:t:pt+3q:1]\qquad p=\frac{1\pm \sqrt{-3}}{2},\quad (t\in \CC).
\end{aligned}
\]

\subsubsection*{Type $\infty$ components}

All components of type $\infty$ are given by \eqref{eq:example_2_1}.

\section{Examples for $n=6$ and $g=1$}\label{sec:more_examples}

\subsection{$(m,n)=(2,6)$ and $g=1$ case}

When $(m,n)=(2,6)$ and $g=1$, the projectivized system is given by
\[
\begin{bmatrix}
a_0b_1\\
a_0b_2+a_1b_1\\
a_0b_3+a_1b_2+a_2b_1\\
a_0b_4+a_1b_3+a_2b_2\\
a_0b_5+a_1b_4+a_2b_3\\
a_1b_5+a_2b_4\\
a_2b_5
\end{bmatrix}
=b_0
\begin{bmatrix}
qa_2\\
a_2-pa_1\\
-pa_2\\0\\0\\0\\0
\end{bmatrix}
-b_0
\begin{bmatrix}
0\\0\\0\\0\\qa_0\\qa_1-a_0\\-pa_0
\end{bmatrix}.
\]

A complete list of solutions are given below.
(We omit the $b$-coordinates from the following list, for brevity.)

\subsubsection*{Type ES components}

Let $\omega$ be a primitive 6th root of the unity.
Define $D_j:=p+q\omega^j$ for $j=0,1,2,3,4,5$.
Solving the system of equations, we have the following components of type ES.

\[
\begin{aligned}
[a_0:a_1:a_2]=&[1:2:1],[1:\frac{-1\pm \sqrt{5}}{2}:1],
[1:\zeta_j:\omega^j],\quad (j=1,2,3,4,5),
\end{aligned}
\]
where $\zeta_j$ is a arbitrary root of the polynomials 
\[
\begin{cases}
D_j\zeta^2-\omega^j\zeta-D_j\omega^j, & (j=1,3,5),\\
D_j\zeta^3-\omega^j \zeta^2-3D_j\omega^j\zeta+2\omega^{2j},& (j=2,4)
\end{cases}
\]
Then, the total number of solutions is $1+2+2+3+2+3+2=15=\binom{6}{2}$.

Let $\epsilon:=p-q$.
At $p=q\to \frac{1}{2}$ limit, these solutions are expressed follows:
\[
\begin{aligned}
[a_0:a_1:a_2]=&[1:2:1],[1:\frac{-1\pm \sqrt{5}}{2}:1],\\
&[1:1+\omega^j:\omega^j],[1:-\frac{1+\omega^j}{3}:\omega^j],\quad (j=1,5),\\
&[1:1+\omega^{j}:\omega^j],[1:(1+\omega^{j})\frac{1\pm \sqrt{17}}{2}:\omega^j],\quad (j=2,4),\\
&\lim_{\epsilon\to 0}[1:\frac{-1\pm \sqrt{1-4\epsilon^2}}{2\epsilon }:-1],\quad (j=3).
\end{aligned}
\]

\subsubsection*{Type X components}

$S$ does not contain irreducible components of type X.


\subsubsection*{Type $\infty$ components}

\[
[a_0:a_1:a_2]=[0:1:0].
\]

Figure \ref{fig:example_24} shows a complete list of solutions, where the eigenvalues are calculated by 
\[
\lambda=\sum_{i=1}^2(qz_1-1+pz_i^{-1})=-2+q\frac{a_1}{a_0}+p\frac{a_1}{a_2}. 
\]

\begin{figure}[htbp]
\begin{center}
\begin{tabular}{c||c|c|c}
& $\{z_1,z_2,z_3\}$& Eigenvalue & XXX model\\\hline
A & $\{1,1\}$ & $0$ & contains $2$ $\infty$'s\\
$\mathrm{B}^{\pm}$ & $\{b^\pm_{1},b^\pm_{2}\}$ & $\frac{-5\pm \sqrt{5}}{2}$ & finite\\
$\mathrm{C}_j$ $(j=1,5)$ & $\{1,\omega^j\}$& $-\frac{1}{2}$ & contains $1$ $\infty$\\
$\mathrm{C}_j$ $(j=2,4)$ & $\{1,\omega^j\}$& $-\frac{3}{2}$ & contains $1$ $\infty$\\
$\mathrm{D}_j$ $(j=1,5)$ & $\{d_{j1},d_{j2}\}$&$-\frac{5}{2}$ & finite\\
$\mathrm{E}^\pm_j$ $(j=2,4)$ & $\{e^\pm_{j1},e^\pm_{j2}\}$&
$\frac{-7\pm \sqrt{17}}{4}$
&finite
\\
$\mathrm{F}$ & $\{1,-1\}$&$-2$ & contains $1$ $\infty$\\
$\mathrm{G}$ & $\{0,\infty\}$&$\frac{-3+\sqrt{1-4\epsilon^2}}{2}\to -1$ &
\begin{tabular}{c}
finite (physical singular)
\end{tabular}
\\
\end{tabular}
\end{center}
\caption{A complete list of solutions to the ASEP-Bethe equations with $(m,n)=(2,6)$, $p=q=\frac{1}{2}$, and $g=1$.
Here, $b_i^{\pm}$ is a root of $b^2-\frac{-1\pm\sqrt{5}}{2}b+1$,
$d_{ji}$ is a root of $d_j^2+\frac{1+\omega^j}{3}d_j+\omega^j$, 
and $e^{\pm}_{ji}$ is a root of $e^2-(1+\omega^j)\frac{1\pm\sqrt{17}}{2}e+\omega^j$.
Solutions $\{z_1,z_2\}$ containing $1$ correspond to infinite solutions to the XXX-Bethe equations.  
The number of finite solutions to the XXX-Bethe equations is $9$.}
\label{fig:example_24}
\end{figure}

\subsection{$(m,n)=(3,6)$ and $g=1$ case}

When $(m,n)=(3,6)$ and $g=1$, the projectivized problem is given by
\begin{equation*}
\begin{bmatrix}
a_0b_1\\
a_0b_2+a_1b_1\\
a_0b_3+a_1b_2+a_2b_1\\
a_0b_4+a_1b_3+a_2b_2+a_3b_1\\
a_0b_5+a_1b_4+a_2b_3+a_3b_2\\
a_1b_5+a_2b_4+a_3b_3\\
a_2b_5+a_3b_4\\
a_3b_5\\
\end{bmatrix}
=b_0
\begin{bmatrix}
q^2a_3\\
2 q a_3-p q a_2\\
p^2 a_1 - p a_2 + (1 - p q) a_3\\
p^2 a_2-p a_3\\
p^2 a_3\\
0\\0\\0
\end{bmatrix}
+
b_0
\begin{bmatrix}
0\\0\\0\\
q^2a_0\\
q^2a_1-qa_0\\
q^2a_2-qa_1+(1-pq)a_0\\
2pa_0-pqa_1\\
p^2a_0
\end{bmatrix}.
\end{equation*}

\subsubsection*{Type ES components}

For a primitive 6th root of the unity $\omega$, let 
\[
\begin{gathered}
\Delta_j:=p^3\omega^{-j}+q^3\omega^j+pq\zeta_j,\quad
\sigma_j:=(\omega^j+\omega^{-j})^{-1},\\
A_j=q^2\omega^j(\zeta_j+1)+p\sigma_j\left\{\omega^j(\zeta_j+pq)+
\omega^{-j}(\zeta_j^2-pq)
\right\},\\
B_j=p^2(\zeta_j+1)+q\sigma_j\left\{\omega^{2j}(\zeta_j^2-pq)+
\zeta_j+pq
\right\}
\end{gathered}
\]
where $\zeta_j$ $(j=0,1,2,3,4,5)$ is a arbitrary root of
\begin{equation}\label{eq:example_zeta}
\begin{cases}
(\zeta-2)(\zeta^3+2\zeta^2+(3-12pq)\zeta + 4pq-2), & (j=0),\\
\zeta^3-\zeta^2-3pq\zeta+4pq-1, & (j=1,5),\\
\zeta^3-\zeta^2-3pq\zeta-2pq+1, & (j=2,4),\\
\zeta^4-(1+12pq)\zeta^2+(8-20pq)\zeta+16pq-4, & (j=3).
\end{cases}
\end{equation}
Then, the following is a complete list of solutions of type E:
\[
[a_0:a_1:a_2:a_3]=[1:\frac{A_j}{\Delta_j}:
\frac{B_j}{\Delta_j}:\omega^j
]
\]
The number of solutions is $20=\binom{6}{3}$.

Let $\epsilon:=p-q$. 
At $p=q\to \frac{1}{2}$ limit, we have
\[
\zeta_j=\begin{cases}
2,\frac{-1\pm\sqrt{5}}{2},-1, & (j=0),\\
0,\frac{3}{2},-\frac{1}{2},
& (j=1,5),\\
\frac{1\pm\sqrt{17}}{4},\frac{1}{2},
& (j=2,4),\\
0,\frac{-1\pm \sqrt{13}}{2},1,
& (j=3),
\end{cases}\qquad
\Delta_j=
\frac{1}{4}\times
\begin{cases}
\zeta_j+1, & (j=0),\\
\zeta_j+\frac{1}{2}, & (j=1,5),\\
\zeta_j-\frac{1}{2}, & (j=2,4),\\
\zeta_j-1, & (j=3).
\end{cases}
\]
If $\Delta_j\neq 0$, we have the polynomial expressions for $a_1/a_0$ and $a_2/a_0$:
\[
\frac{a_1}{a_0}=
\begin{cases}
\zeta_j+1, & (j=0),\\
2(1-\omega^j)\zeta_j-1+4\omega^j, & (j=1,5),\\
2(1+\omega^{j})\zeta_j+1, & (j=2,4),\\
\zeta_j+1, & (j=3).
\end{cases}\qquad
\frac{a_2}{a_0}=
\begin{cases}
\zeta_j+1, & (j=0),\\
-2(1-\omega^j)\zeta_j+4-\omega^j, & (j=1,5),\\
2(1+\omega^j)\zeta_j+\omega^j, & (j=2,4),\\
-1-\zeta_j, & (j=3).
\end{cases}
\]
If $\Delta_j= 0$, we need a perturbation:
\begin{itemize}
\item Around $\zeta_0=-1$, we have the expansion $\zeta_0=-1-4\epsilon^2-28\epsilon^4+\cdots$, $\Delta_0=-6\epsilon^4-86\epsilon^6+\cdots$, $A_0=3\epsilon^3+3\epsilon^4+\cdots$, and $B_0=-3\epsilon^3+3\epsilon^4+\cdots$.
Then, it follows that $a_1/a_0= O(\epsilon^{-1})$ and $a_2/a_0=O(\epsilon^{-1})$.
\item Around $\zeta_i=-\frac{1}{2}$ $(j=1,5)$, we have the expansion $\zeta_j=-\frac{1}{2}+\frac{11}{8}\epsilon^2+\cdots$, $\Delta_j=\frac{3}{8} (\omega^{-j}-\omega^j)\epsilon+\cdots$, $A_j=-\frac{3}{8}\omega^j\epsilon+\cdots$, and $B_j=\frac{3}{8}\epsilon+\cdots$.
Then, it follows that $a_1/a_0\to \frac{\omega^{2j}}{\omega^{2j}-1}$ and $a_2/a_0\to \frac{-\omega^{j}}{\omega^{2j}-1}$.
\item Around $\zeta_i=\frac{1}{2}$ $(j=2,4)$, we have the expansion $\zeta_j=\frac{1}{2}+\frac{7}{8}\epsilon^2+\cdots$, $\Delta_j=\frac{3}{8} (\omega^{-j}-\omega^j)\epsilon+\cdots$, $A_j=-\frac{9}{8}\omega^j\epsilon+\cdots$, and $B_j=\frac{9}{8}\epsilon+\cdots$.
Then, it follows that $a_1/a_0\to \frac{3\omega^{2j}}{\omega^{2j}-1}$ and $a_2/a_0\to \frac{-3\omega^{j}}{\omega^{2j}-1}$.
\item Around $\zeta_3=1$, we have the expansion $\zeta_3=1+4\epsilon^2+76\epsilon^4+\cdots$, $\Delta_3=18\epsilon^4+570\epsilon^6+\cdots$, $A_3=\frac{3}{2}\epsilon+\frac{3}{2}\epsilon^2+5\epsilon^3+41\epsilon^4+\cdots$, and $B_3=\frac{3}{2}\epsilon-\frac{3}{2}\epsilon^2+5\epsilon^3-41\epsilon^4+\cdots$.
Then, it follows that $a_1/a_0= O(\epsilon^{-3})$ and $a_2/a_0=O(\epsilon^{-3})$.
\end{itemize}


\subsubsection*{Type X component}

\[
[a_0:a_1:a_2:a_3]=[1:q^{-1}t:p^{-1}(3-t):1],\quad
p/q=\omega^j,\quad (t\in \CC),\qquad (j=2,4).
\]

\subsubsection*{Type $\infty$ component}

\[
[a_0:a_1:a_2:a_3]=[0:s:t:0],
\]
where $pq(t^6-s^6)-pst^5+qs^5t=0$.

Figure \ref{fig:example_36} shows a complete list of solutions, where the eigenvalues are calculated by 
\[
\lambda=\sum_{i=1}^3(qz_i-1+pz_i^{-1})=-3+q\frac{a_1}{a_0}+p\frac{a_2}{a_3}. 
\]

\begin{figure}[htbp]
\begin{center}
\begin{tabular}{c||c|c|c|c}
& $\zeta_j$& $\{z_1,z_2\}$& Eigenvalue & XXX model\\\hline
A $(j=0)$ & $2$ & $\{1,1,1\}$ & $0$ & contains $3$ $\infty$'s\\
$\mathrm{B}^{\pm}$ &$\frac{-1\pm\sqrt{5}}{2}$& $\{1,b^\pm_{1},b^\pm_{2}\}$ & $\frac{-5\pm \sqrt{5}}{2}$ & contains $1$ $\infty$\\
$\mathrm{C}_j$ $(j=1,5)$ &$0$& $\{c_{j1},c_{j2},c_{j3}\}$& $-2$ & finite\\
$\mathrm{D}_j$ $(j=1,5)$ &$\frac{3}{2}$& $\{1,1,\omega^j\}$& $-\frac{1}{2}$ & contains $2$ $\infty$'s\\
$\mathrm{D}_j$ $(j=2,4)$ &$\frac{1}{2}$& $\{1,1,\omega^j\}$& $-\frac{3}{2}$ & contains $2$ $\infty$'s\\
$\mathrm{E}^\pm_j$ $(j=2,4)$ & $\frac{1\pm\sqrt{17}}{4}$& $\{1,e_{j1},e_{j2}\}$&$\frac{-7\pm \sqrt{17}}{4}$ & contains $1$ $\infty$\\
F &$0$& $\{1,1,-1\}$&
$-2$
& contains $2$ $\infty$'s
\\
$\mathrm{G}^\pm$ &$\frac{-1\pm\sqrt{13}}{2}$& $\{-1,g_1,g_2\}$&$\frac{-5\pm\sqrt{13}}{2}$ & finite\\
$\mathrm{H}$ &$-1$& $\{0,-1,\infty\}$&$-3$ & finite (physical singular)\\
$\mathrm{I}_{j}$ $(j=1,5)$ & $-\frac{1}{2}$ & $\{1,i_{j1},i_{j2}\}$ & $-\frac{5}{2}$ & contains $1$ $\infty$\\
$\mathrm{J}$& $1$ & $\{0,1,\infty\}$ & $-1$ & contains $1$ $\infty$ (physical singular\footnote{})\\
\end{tabular}
\end{center}
\caption{A complete list of solutions to the ASEP-Bethe equations, arising from type ES components, for $(m,n)=(3,6)$, $p=q=\frac{1}{2}$, and $g=1$.
Here, $b_i^\pm,c_{ji},e_{ji},g_i,i_{ji}$ are roots of polynomials associated to $\zeta_j$ given in \eqref{eq:example_zeta}.
The number of finite solutions to the XXX-Bethe equations is
$5$.\\[.5em]
{\small ${}^4$ According to the original definition of the term, this is not a ``physical singular solution'' because Solution J corresponds to the XXX-Bethe root $\{u_1,u_2,u_3\}=\{-\frac{i}{2},\infty,\frac{i}{2}\}$, which contains $\infty$.
However, by its nature, it should be appropriate to say that ``Solution J is a physical singular solution with one $\infty$.'' 
}}
\label{fig:example_36}
\end{figure}


\newpage

\begin{thebibliography}{10}

\bibitem{alba2013bethe}
Vincenzo Alba, Kush Saha, and Masudul Haque.
\newblock Bethe ansatz description of edge-localization in the open-boundary
  {XXZ} spin chain.
\newblock {\em Journal of Statistical Mechanics: Theory and Experiment},
  2013(10):P10018, 2013.

\bibitem{avdeev1986exceptional}
LV~Avdeev and Artem~Alekseevich Vladimirov.
\newblock Exceptional solutions to the {B}ethe ansatz equations.
\newblock {\em Teoreticheskaya i Matematicheskaya Fizika}, 69(2):163--174,
  1986.

\bibitem{brattain2015completeness}
Eric Brattain, Norman Do, and Axel Saenz.
\newblock The completeness of the {B}ethe ansatz for the periodic {ASEP}.
\newblock {\em arXiv preprint arXiv:1511.03762}, 2015.

\bibitem{crampe2018equivalent}
Nicolas Crampe and Rafael~I Nepomechie.
\newblock Equivalent {TQ} relations and exact results for the open {TASEP}.
\newblock {\em Journal of Statistical Mechanics: Theory and Experiment},
  2018(10):103105, 2018.

\bibitem{crampe2010eigenvectors}
Nicolas Cramp{\'e}, Eric Ragoucy, and Damien Simon.
\newblock Eigenvectors of open {XXZ} and {ASEP} models for a class of
  non-diagonal boundary conditions.
\newblock {\em Journal of Statistical Mechanics: Theory and Experiment},
  2010(11):P11038, 2010.

\bibitem{de2006exact}
Jan De~Gier and Fabian~HL Essler.
\newblock Exact spectral gaps of the asymmetric exclusion process with open
  boundaries.
\newblock {\em Journal of Statistical Mechanics: Theory and Experiment},
  2006(12):P12011--P12011, 2006.

\bibitem{deguchi2001non}
Tetsuo Deguchi.
\newblock Non-regular eigenstate of the {XXX} model as some limit of the bethe
  state.
\newblock {\em Journal of Physics A: Mathematical and General},
  34(46):9755--9775, 2001.

\bibitem{fantechi2005fundamental}
Barbara Fantechi.
\newblock {\em Fundamental algebraic geometry: Grothendieck's FGA explained}.
\newblock Number 123. American Mathematical Soc., 2005.

\bibitem{fulton2013intersection}
William Fulton.
\newblock {\em Intersection Theory}.
\newblock Ergebnisse der Mathematik und ihrer Grenzgebiete. 3. Folge / A Series
  of Modern Surveys in Mathematics. Springer Berlin Heidelberg, 2013.

\bibitem{golinelli2006asymmetric}
Olivier Golinelli and Kirone Mallick.
\newblock The asymmetric simple exclusion process: an integrable model for
  non-equilibrium statistical mechanics.
\newblock {\em Journal of Physics A: Mathematical and General},
  39(41):12679--12705, 2006.

\bibitem{hartshorne2013algebraic}
Robin Hartshorne.
\newblock {\em Algebraic geometry}.
\newblock Springer Science \& Business Media, 2013.

\bibitem{iwao2025bethe}
Shinsuke Iwao and Kohei Motegi.
\newblock Bethe roots for periodic {TASEP} and algebraic curve.
\newblock {\em arXiv preprint arXiv:2504.19690}, 2025.

\bibitem{jiang2018algebraic}
Yunfeng Jiang and Yang Zhang.
\newblock Algebraic geometry and {B}ethe ansatz. {P}art {I}. the quotient ring
  for {BAE}.
\newblock {\em Journal of High Energy Physics}, 2018(3):1--40, 2018.

\bibitem{kirillov2014singular}
Anatol~N Kirillov and Reiho Sakamoto.
\newblock Singular solutions to the {B}ethe ansatz equations and rigged
  configurations.
\newblock {\em Journal of Physics A: Mathematical and Theoretical},
  47(20):205207, 2014.

\bibitem{KirillovSakamoto2015}
Anatol~N. Kirillov and Reiho Sakamoto.
\newblock Some remarks on {N}epomechie--{W}ang eigenstates for spin 1/2 {XXX}
  model.
\newblock {\em Moscow Mathematical Journal}, 15(2):337--352, 2015.

\bibitem{levkovich2016bethe}
Fedor Levkovich-Maslyuk.
\newblock The {B}ethe ansatz.
\newblock {\em Journal of Physics A: Mathematical and Theoretical},
  49(32):323004, 2016.

\bibitem{matsumura1989commutative}
Hideyuki Matsumura.
\newblock {\em Commutative ring theory}.
\newblock Number~8. Cambridge university press, 1989.

\bibitem{nepomechie2013algebraic}
Rafael~I Nepomechie and Chunguang Wang.
\newblock Algebraic {B}ethe ansatz for singular solutions.
\newblock {\em Journal of Physics A: Mathematical and Theoretical},
  46(32):325002, 2013.

\bibitem{nepomechie2014twisting}
Rafael~I Nepomechie and Chunguang Wang.
\newblock Twisting singular solutions of {B}ethe's equations.
\newblock {\em Journal of Physics A: Mathematical and Theoretical},
  47(50):505004, 2014.

\bibitem{noh2000origin}
Jae~Dong Noh, Deok-Sun Lee, and Doochul Kim.
\newblock Origin of the singular {B}ethe ansatz solutions for the {H}eisenberg
  {XXZ} spin chain.
\newblock {\em Physica A: Statistical Mechanics and its Applications},
  287(1-2):167--176, 2000.

\bibitem{prolhac2020riemann}
Sylvain Prolhac.
\newblock Riemann surface for {TASEP} with periodic boundaries.
\newblock {\em Journal of Physics A: Mathematical and Theoretical},
  53(44):445003, 2020.

\bibitem{prolhac2020riemannKPZ}
Sylvain Prolhac.
\newblock Riemann surfaces for {KPZ} with periodic boundaries.
\newblock {\em SciPost Physics}, 8(1):008, 2020.

\bibitem{prolhac2021riemann}
Sylvain Prolhac.
\newblock From the {R}iemann surface of {TASEP} to {ASEP}.
\newblock {\em Journal of Physics A: Mathematical and Theoretical},
  54(39):395002, 2021.

\end{thebibliography}

\end{document}